\documentclass[10pt]{amsart}
\usepackage{diagrams}
\usepackage{amssymb}
\usepackage{amsmath}

\newtheorem{proposition}{Proposition}[section]
\newtheorem{lemma}[proposition]{Lemma}
\newtheorem{corollary}[proposition]{Corollary}
\newtheorem{theorem}[proposition]{Theorem}

\theoremstyle{definition}

\newtheorem{example}[proposition]{Example}

\theoremstyle{remark}

\newcommand{\thlabel}[1]{\label{th:#1}}
\newcommand{\thref}[1]{Theorem~\ref{th:#1}}
\newcommand{\selabel}[1]{\label{se:#1}}

\newcommand{\lelabel}[1]{\label{le:#1}}
\newcommand{\leref}[1]{Lemma~\ref{le:#1}}
\newcommand{\prlabel}[1]{\label{pr:#1}}
\newcommand{\prref}[1]{Proposition~\ref{pr:#1}}
\newcommand{\colabel}[1]{\label{co:#1}}
\newcommand{\coref}[1]{Corollary~\ref{co:#1}}

\newcommand{\exlabel}[1]{\label{ex:#1}}
\newcommand{\exref}[1]{Example~\ref{ex:#1}}

\def\hp{\hspace{.5cm}}
\def\<{\leq}
\def\>{\geq}

\def\b{\beta}

\def\D{\Delta}

\def\e{\varepsilon}

\def\s{\sigma}
\def\S{\sum}
\def\u{\underline}
\def\us{\u{\s}}

\def\ot{\otimes}
\def\ra{\rightarrow}
\def\rhu{\rightharpoonup}

\def\BQ{\mathrm{BQ}}
\def\BM{\mathrm{BM}}

\def\no{n_{(0)}}
\def\al{a_{(1)}}
\def\bl{b_{(1)}}
\def\az{a_{(2)}}
\def\bz{b_{(2)}}

\def\nl{n_{(1)}}

\def\BQ{\mathrm{BQ}}
\def\BM{\mathrm{BM}}
\def\BC{\mathrm{BC}}
\def\BW{\mathrm{BW}}
\def\YDH{_{H}\mathcal{YD}^H}
\def\YDHs{{_{H^{\s}}\mathcal{YD}}^{H^{\s}}}
\def\thetaH{{H_{\theta}}}
\def\HR{\mathcal{H}_R}
\def\Gal{\mathrm{Gal}}
\def\Br{\mathrm{Br}}
\date{}

\begin{document}
\title{Cocycle Deformations and Brauer Group Isomorphisms}
\author{Hui-xiang Chen}
\address{Department of Mathematics, Yangzhou University,
Yangzhou 225002, China\\
\& School of Mathematics, Statistics and Computer Science
Victoria University of Wellington,
P.O. Box 600, Wellington, New Zealand}
\email{hchen@mcs.vuw.ac.nz, yzchenhx@yahoo.com}
\author{Yinhuo Zhang}
\address{School of Mathematics, Statistics and Computer Science,
Victoria  University of Wellington, PO Box 600,
Wellington, New Zealand}
\email{yinhuo.zhang@vuw.ac.nz}
\begin{abstract}
Let $H$ be a Hopf algebra over a commutative ring $k$ with unity and $\sigma:H\otimes H\longrightarrow k$ be a cocycle on $H$. In this paper, we show that the Yetter-Drinfeld module category of the cocycle deformation Hopf algebra $H^{\sigma}$ is equivalent to the Yetter-Drinfeld module category of $H$. As a result of the equivalence, the ``quantum Brauer" group BQ$(k,H)$ is isomorphic to BQ$(k,H^{\sigma})$. Moreover, the group $\Gal(\HR)$ constructed in \cite{Z} is studied under a cocycle deformation.
\end{abstract}
\maketitle

\def\r{\longrightarrow}
\def\M{\mathcal M}
\section*{\bf Introduction}
Let $k$ be a commutative ring with unity and let $H$ be a Hopf algebra over $k$. In \cite{Doi}, Doi introduced a cocycle twisted Hopf algebra $H^{\sigma}$, called a cocycle deformation of $H$, for a 2-cocycle $\sigma:~H\otimes H\r k$. It is shown in \cite{DT} that Drinfeld's quantum double $D(H)$ is a cocycle deformation of the tensor product Hopf algebra $H^{*cop}\otimes H$. If $(H,R)$ is a coquasitriangular Hopf algebra, then $(H^{\sigma}, R^{\sigma})$ is a coquasitriangular Hopf algebra as well (see \cite[p61]{Maj}), where $R^{\sigma}=(\sigma\tau) R\sigma^{-1}$ and $\tau$ is the flip map. It is well known that the braided monoidal category $\M^H_R$ of right $H$-comodules with the braiding induced by the coquasitriangular structure $R$ is equivalent to the braided monoidal category $\M^{H^{\sigma}}_{R^{\sigma}}$ of right $H^{\sigma}$-comodule category with the braiding induced by $R^{\sigma}$. Since the Brauer group Br$(\mathcal C)$ is a group invariant of a braided monoidal category $\mathcal{C}$ \cite{VZ0}, we have a group isomorphism from the equivariant Brauer group  BC$(k,H,R)$ to the equivariant Brauer group BC$(k,H^{\sigma}, R^{\sigma})$ \cite{Car}.
The significance of this isomorphism was demonstrated via Sweedler's Hopf algebra $H_4$ which has a family of coquasitrianglar structures $R_t$ indexed by the ground ring $k$. Since each $(H_4, R_t)$ is a cocycle deformation of $(H_4, R_0)$, we have group isomorphisms BC$(k,H_4,R_t)\cong \mathrm{BC}(k,H_4,R_0)$, for all $t\in k$, where the later is known as the direct product group of the Brauer-Wall group BW$(k)$ and the additive group $k^+$ \cite{VZ1}. In this paper, we will show that the braided monoidal category equivalence between $\M^H_R$ and $\M^{H^{\sigma}}_{R^{\sigma}}$ is carried out by an equivalence braided monoidal functor from the Yetter-Drinfeld module category $_{H}{\mathcal YD}^H$ to the Yetter-Drinfeld module category $_{H^{\sigma}}{\mathcal YD}^{H^{\sigma}}$. As a consequence, we obtain a Brauer group isomorphism between BQ$(k,H)$ and BQ$(k,H^{\sigma})$ for any cocycle deformation $H^{\sigma}$ of a Hopf algebra $H$ with a bijective antipode. When $(H,R)$ is a coquasistriangular Hopf algebra, then the group isomorphism restricts to the group isomorphism between the equivariant Brauer groups BC$(k,H,R)$ and BC$(k,H^{\sigma}, R^{\sigma})$ obtained in \cite{Car}. The equivalence of Yetter-Drinfeld module categories under cocycle deformation is of self-duality in the sense that the Yetter-Drinfeld module category $_{H}{\mathcal YD}^H$ is not only stable under cocycle deformation but also stable under dual-cocycle deformation. In other words, if $\theta\in H\otimes H$ is a  dual-cocycle and $\thetaH $ is the dual-cocycle deformation of $H$, then the braided monoidal category $_{H}{\mathcal YD}^H$ is equivalent to the braided monoidal category $_{\thetaH }{\mathcal YD}^{\thetaH }$.

In the first section, we recall some basic definitions about the Yetter-Drinfeld module category of a Hopf algebra and the (equivariant) Brauer group of a (coquasitriangular) Hopf algebra. In Section 2, we show that for a Hopf algebra $H$ there exists a natural equivalence braided monoidal functor $\underline{\sigma}$ from the Yetter-Drinfeld module category $_{H}{\mathcal YD}^H$  to the Yetter-Drinfeld module category $_{H^{\sigma}}{\mathcal YD}^{H^{\sigma}}$ of the cocycle deformation $H^{\sigma}$ with respect to a  cocycle $\sigma:H\otimes H\r k$ (see Theorem 2.3); and there exists a natural equivalence braided monoidal functor from  $_{H}{\mathcal YD}^H$ to the Yetter-Drinfeld module category $_{\thetaH }{\mathcal YD}^{\thetaH }$ of a dual-cocycle deformation $\thetaH $ with respect to the  dual cocycle $\theta\in H\otimes H$. When $(H,R)$ is a coquasitriangular Hopf algebra, then the equivalence  functor $\underline{\sigma}$ restricts to an equivalence functor from the full braided monoidal subcategory $\M^H_R$ of right $H$-comodules to the full braided monoidal subcategory $\M^{H^{\sigma}}_{R^{\sigma}}$. It is well-known that the Brauer group BQ$(k,H)$ is the Brauer group of the braided monoidal category $_{H}{\mathcal YD}^H$. Thus we obtain a group isomorphism from BQ$(k,H)$ to BQ$(k,H^{\sigma})$. When $\sigma$ is a lazy cocycle, the Hopf algebra $H^\sigma$ is equal to $H$ and $\sigma$ induces an automorphism of the Brauer group $\BQ(k,H)$. We obtain that the second lazy cohomology group $\mathrm{H}^2_L(H)$ acts on $\BQ(k,H)$ by automorphisms (see Corollary 2.7).

Section 3 devotes to the study of the group Gal$(\HR)$ of bigalois objects constructed in \cite{Z} under cocycle deformation. If $(H,R)$ is a finite (faithfully projective) coquasitriangular Hopf algebra, we have a generalized cotensor product over the braided Hopf algebra $\HR$ defined in the Yetter-Drinfeld module category $\YDH$ (see \cite{Z} for more detail). The generalized cotensor product induces a group structure on the set Gal$(\HR)$ of quantum commutative $\HR$-bigalois objects in the category $\YDH$. The group Gal$(\HR)$ measures the equivariant Brauer group BC$(k,H,R)$ via an exact group sequence:
\begin{equation}\label{exactsqn}
1\r \Br(k)\r \BC(k,H,R)\r \Gal(\HR).
\end{equation}
Now let $\sigma$ be a  cocycle on $H$ and let $(H^{\sigma}, R^{\sigma})$ be the cocycle deformation of $(H,R)$ with respect to $\sigma$.
We show that the equivalence functor $\underline{\sigma}$ defined in Section 2 commutes with the generalized cotensor product, and induces a group isomorphism from Gal$(\HR)$ to Gal$(\mathcal{H}^{\sigma}_{R^{\sigma}})$. Moreover, the exact sequence (\ref{exactsqn}) is stable under cocycle deformation. That is,  we have the following commutative diagram of exact sequences:
$$\begin{diagram}
1 &\rTo &  \Br(k)&\rTo & \BC(k,H,R)& \rTo^{\tilde{\pi}} & \Gal(\HR)\\
&& \dTo_{=} &  & \dTo_{\us} & & \dTo_{\us} \\
1 &\rTo &  \Br(k)&\rTo & \BC(k,H^{\sigma},R^{\sigma})& \rTo^{\tilde{\pi}} & \Gal(\mathcal{H}^{\s}_{R^{\sigma}}).
\end{diagram}$$

\newpage
\section{\bf Preliminaries}\selabel{1}
\subsection{Yetter-Drinfeld modules}\selabel{1.1}
Throughout, we work over a fixed commutative ring $k$ with unity. Unless
otherwise stated, all algebras,  Hopf algebras and modules are
defined over $k$; all maps are $k$-linear; dim,
$\otimes$ and Hom stand for $\mbox{\rm dim}_k$, $\otimes_k$ and
Hom$_k$, respectively. For the theory of Hopf algebras and quantum groups, we
refer to \cite{Ka, Maj, Mon, Sw}.

Let $H$ be a Hopf algebra with a bijective antipode $S$.
A Yetter-Drinfeld $H$-module (simply, YD $H$-module) $M$ is a crossed $H$-bimodule.
That is, $M$ is at once a left $H$-module and a right $H$-comodule satisfying
the following equivalent compatibility conditions:
\begin{eqnarray}\label{compat}
\Sigma h_{(1)}\cdot m_{(0)}\ot h_{(2)}m_{(1)}
&=& \Sigma (h_{(2)}\cdot m)_{(0)}\ot (h_{(2)}\cdot m)_{(1)}h_{(1)},\\
\Sigma(h\cdot m)_{(0)}\ot(h\cdot m)_{(1)}
&=& \Sigma h_{(2)}\cdot m_{(0)}\ot h_{(3)}m_{(1)}S^{-1}(h_{(1)}),
\end{eqnarray}
where $h\in H$, $m\in M$, and the sigma notations for a comodule and for a comultiplication can be found in the reference book \cite{Sw}. The category $_{H}{\mathcal{YD}}^H$ of Yetter-Drinfeld $H$-modules and their homomorphisms is a braided monoidal category with the braiding given by
$$\Phi: M\otimes N\r N\otimes M, \ \Phi(m\otimes n)=\sum n_{(0)}\otimes n_{(1)}\cdot m$$
where $M,N\in {\YDH}$.

A Hopf algebra over $k$ is called a coquasitriangular (CQT) Hopf algebra if there is an invertible element $R\in (H\ot H)^*$, the convolution
algebra of $H\ot H$, subject to the following conditions:

(CQT1). $R(h\ot 1)=R(1\ot h)=\varepsilon(h)$,\\
(CQT2). $R(g\ot hl)=\sum R(g_{(1)}\ot l) R(g_{(2)}\ot h)$,\\
(CQT3). $R(hl\ot g)=\sum R(h\ot g_{(1)}) R(l\ot g_{(2)})$,\\
(CQT4). $\sum R(g_{(1)}\ot h_{(1)})g_{(2)}h_{(2)} =\sum R(g_{(2)}\ot h_{(2)})h_{(1)}g_{(1)}$,

for all $g, h$ and $l\in H$. Now let $H$ be a CQT Hopf algebra with a CQT structure $R$. Note that (CQT4) is equivalent to one of the following:

(CQT$4'$).  $\S g_{(1)}R(g_{(2)}\ot h)=\S R(g_{(1)}\ot h_{(2)})S(h_{(1)})g_{(2)}h_{(3)},
\hspace{0.1cm}g,h\in H$,\\
(CQT$4''$).  $\S h_{(1)}R(g\ot h_{(2)})=\S R(g_{(2)}\ot h_{(1)})g_{(3)}h_{(2)}S^{-1}(g_{(1)}),
\hspace{0.1cm}g,h\in H$.

If $M$ is a right $H$-(or $H^{op}$)-comodule, the CQT structure $R$ induces a left $H$-module structure on $M$ as follows:
\begin{equation}\label{Raction}
h\vartriangleright_1 m=\sum m_{(0)} R(h\ot m_{(1)})
\end{equation}
for $h\in H$ and $m\in M$, where we use the $\vartriangleright_1$ as we will have the second left $H$-action induced by the CQT structure $R$ in Section 2. The $H$-action (\ref{Raction}) together with the original $H$-coaction makes $M$ into a YD $H$-module. Denote by $\M^H_R$ the category of YD $H$-modules with the left $H$-module structure (\ref{Raction})
coming from the right $H$-comodule structure. It is obvious
that $\M^H_R$ is a full braided monoidal subcategory of $_H{\mathcal{YD}}^H$.

\def\R{\mathcal{R}}
Now let $H$ be a quasitriangular (QT) Hopf algebra, that is, $H$ is a Hopf algebra
with an invertible element $\R=\S \R^{(1)}\ot \R^{(2)}$ in $H\ot H$ satisfying
the following axioms ($r=\R$):

(QT1). \hp  $\sum\D(\R^{(1)})\ot \R^{(2)} =\sum \R^{(1)}\ot r^{(1)}\ot \R^{(2)} r^{(2)},$\\
(QT2). \hp  $\sum\e(\R^{(1)}) \R^{(2)} =\sum \R^{(1)}\e(\R^{(2)})=1$,\\
(QT3). \hp  $\sum  \R^{(1)}\ot \Delta( \R^{(2)})=\sum \R^{(1)} r^{(1)} \ot r^{(2)}\ot  \R^{(2)},$\\
(QT4). \hp  $\D^{cop}(h)\R= \R\D(h),$

where $\Delta^{cop}= \tau\Delta$ is the comultiplication of the Hopf
algebra $H^{cop}$ and $\tau$ is the flip map.

If $A$ is a left $H$-module (algebra), then $A$ is
simultaneously a YD $H$-module (algebra) with the right $H^{op}$-comodule structure
given by
\begin{equation}\label{Rcomodule}
A\ra A\ot H^{op},\hspace{0.2cm} a\mapsto \S (\R^{(2)}\cdot a)\ot \R^{(1)},\hspace{0.2cm}a\in A.
\end{equation}

\subsection{$H$-Azumaya algebras}\selabel{1.2}
A Yetter-Drinfeld $H$-module algebra (simply, YD $H$-module algebra) is a YD $H$-module $A$ such that $A$ is a left $H$-module algebra and a right $H^{\rm op}$-comodule algebra.
For the details of $H$-(co)module algebras we refer to \cite{Ka, Maj, Mon, Sw}.

Let $A$ and $B$ be two YD $H$-module algebras. We may define a {\it
braided product}, still denoted $\#$, on the YD $H$-module $A\ot B$:
$$(a\#b)(c\#d)=\sum ac_{(0)} \# (c_{(1)}\cdot b)d$$
for $a, c\in A$ and $b, d\in B$. The braided product $\#$ makes $A\# B$
a left $H$-module algebra and a right $H^{op}$-comodule algebra so that
$A\# B$ is a YD $H$-module algebra.

Now let $A$ be a YD $H$-module algebra. The {\it $H$-opposite algebra
$\overline{A}$} of $A$ is the YD $H$-module algebra defined as
follows: $\overline{A}$ equals $A$ as a YD $H$-module, but
with multiplication given by the formula
$$\overline{a}\circ\overline{b} = \sum \overline{b_{(0)} (b_{(1)}\cdot a)}$$
for all $\overline{a}, \overline{b} \in \overline{A}$.

\def\End{\mathrm{End}}
\def\ov{\overline}
\def\ydmh{Yetter-Drinfeld $H$-module algebra}
In \cite{CVZ1} we defined the Brauer group of a Hopf algebra $H$ by
considering isomorphism classes of $H$-Azumaya algebras.
A \ydmh\ $A$ is said to be  {\it $H$-Azumaya} if
it is finite (i.e., faithfully projective) as a $k$-module and if the
following  two \ydmh\ maps are isomorphisms:
$$\begin{array}{lllll}
F:& A\# \ov{A} & \r & \End(A), & F(a\# \ov{b})(x)=\sum ax_{(0)} (x_{(1)}\cdot b),\\
G:& \ov{A}\# A & \r & \End(A)^{op}, & G(\ov{a}\# b)(x)=\sum a_{(0)} (a_{(1)}\cdot x)b,
\end{array}$$
For a finite YD $H$-module $M$, the endomorphism algebra $\End_k(M)$
is a \ydmh\ with the $H$-structures given by
$$\begin{array}{l}
(h\cdot f)(m) = \sum h_{(1)}\cdot f(S(h_{(2)})\cdot m),\\
\sum f_{(0)}(m)\otimes f_{(1)}  = \sum f(m_{(0)})_{(0)} \ot
S^{-1}(m_{(1)})f(m_{(0)})_{(1)}
\end{array}$$
for $f\in \End(M)$ and $m\in M$.
The elementary $H$-Azumaya algebra $\End(M)^{op}$ has the different
 $H$-structures from those of $\End(M)$ (see \cite{CVZ1} for the details).

Two $H$-Azumaya algebras $A$ and $B$ are Brauer equivalent (denoted
 $A\sim B$) if there exist two finite YD $H$-modules
$M$ and $N$ such that $A\#\End(M) \cong B\#\End(N)$. Note that $A\sim B$ if
and only if $A$ is $H$-Morita equivalent to $B$ (see \cite[Th.2.10]{CVZ1}).
The relation $\sim$ is an equivalence relation on the set $B(k,H)$ of isomorphism classes of $H$-Azumaya algebras and the quotient set of $B(k,H)$ modulo $\sim$ is a group, called the Brauer group of the Hopf algebra $H$, denoted $\BQ(k,H)$.

If $(H,R)$ is a CQT Hopf algebra, then BQ$(k,H)$ contains a subgroup consisting of elements  represented by $H$-Azumaya algebras whose $H$-module structures are induced by the CQT structure $R$ via (\ref{Raction}). Such an $H$-Azumaya algebra is called an $R$-Azumaya algebra. The subgroup, denoted BC$(k,H,R)$, is called the {\it equivariant Brauer group} of the CQT Hopf algebra $(H,R)$ \cite{CVZ2}.

\section{\bf The Yetter-Drinfeld module category of a cocycle deformation}\selabel{2}

Let $H$ be a Hopf algebra. Recall that a 2-cocycle on $H$ is a convolution invertible $k$-linear map $\s: H\ot H\ra k$ satisfying:
\begin{equation}\label{lcocycle}
\S\s(g_{(1)}\ot h_{(1)})\s(g_{(2)}h_{(2)}\ot l)=
\S\s(h_{(1)}\ot l_{(1)})\s(g\ot h_{(2)}l_{(2)}),
\end{equation}
or equivalently,
\begin{equation}\label{+-}
\S\s(g_{(1)}h_{(1)}\ot l_{(1)})\s^{-1}(g_{(2)}\ot h_{(2)}l_{(2)})
=\S\s^{-1}(g\ot h_{(1)})\s(h_{(2)}\ot l),
\end{equation}
and $\s(h\ot 1)=\s(1\ot h)=\e(h)1$, for all $g$, $h$, $l\in H$.
In general, the inverse $\s^{-1}$ of  a  2-cocycle $\s$ on $H$ is not necessarily  a 2-cocycle on $H$. But $\s^{-1}$ satisfies:
\begin{equation}\label{rcocycle}
\S\s^{-1}(g_{(1)}h_{(1)}\ot l)\s^{-1}(g_{(2)}\ot h_{(2)})=
\S\s^{-1}(g\ot h_{(1)}l_{(1)})\s^{-1}(h_{(2)}\ot l_{(2)})
\end{equation}
for all $g$, $h$, $l\in H$.  We will use the equation (\ref{rcocycle}) later in computations. Furthermore, a 2-cocycle $\s$ on $H$ satisfies the following identity by \cite[Theorem 1.6(a5)]{Doi}:
\begin{equation}\label{+-e}
\S\s(h_{(1)}\ot S(h_{(2)}))\s^{-1}(S(h_{(3)})\ot h_{(4)})=\e(h),\hspace{0.2cm}h\in H.
\end{equation}

In \cite{Doi}, Y. Doi introduced a new Hopf algebra $H^{\s}$, called the $\s$-deformation of $H$ for a 2-cocycle $\sigma$ on $H$. The Hopf algebra $H^{\s}$ is equal to $H$ as a coalgebra. But its multiplication and its antipode are given by
\begin{equation}\label{stime}
h\cdot_{\s}h'=\S\s(h_{(1)}\ot h'_{(1)})h_{(2)}h'_{(2)}\s^{-1}(h_{(3)}\ot h'_{(3)})
\end{equation}
and
$$S^{\s}(h)=\S\s(h_{(1)}\ot S(h_{(2)}))S(h_{(3)})\s^{-1}(S(h_{(4)})\ot h_{(5)}),$$
respectively, where $h,h'\in H$ and $S$ is the antipode of $H$. If $S$ is bijective, then $S^{\s}$ is bijective as well and its inverse $(S^{\s})^{-1}$ is given by
$$(S^{\s})^{-1}(h)=\S\s^{-1}(h_{(5)}\ot S^{-1}(h_{(4)}))S^{-1}(h_{(3)})\s(S^{-1}(h_{(2)})\ot h_{(1)})$$
for all $h\in H^{\s}$ (see \cite{Doi, DT}).

From now on, $H$ will be a Hopf algebra with a bijective antipode $S$
and $\s$ will be a 2-cocycle on $H$ if it is not specified.
We show in this section that the Yetter-Drinfeld module categories
${_H\mathcal{YD}}^H$ and $\YDHs$ are equivalent as braided monoidal categories. To this aim, we first construct a covariant functor from ${_H\mathcal{YD}}^H$ to $\YDHs$.

\begin{lemma}\lelabel{2.1} Let $M$ be a YD $H$-module $M$. Then \\
{\rm (a)}\ $M$ is a left $H^{\s}$-module with the $H^{\s}$-action given by
\begin{equation}\label{sigmaaction}
\begin{array}{rcl}
h\rhu m &=& \S(h_{(2)}\cdot m_{(0)})_{(0)}\s((h_{(2)}\cdot m_{(0)})_{(1)}\ot h_{(1)}) \s^{-1}(h_{(3)}\ot m_{(1)})\\
&=& \S(h_{(3)}\cdot m_{(0)})\s(h_{(4)}m_{(1)}S^{-1}(h_{(2)})\ot h_{(1)})
\s^{-1}(h_{(5)}\ot m_{(2)}),\\
\end{array}
\end{equation}
where $h\in H$ and $m\in M$. \\
{\rm (b)}\ Denote by $\u{\s}(M)$ the left $H^{\s}$-module in {\rm (a)}. If $f:M\r N$ is a {\rm YD} $H$-module map, then $\u{\s}(f)=f:\us(M)\r \us(N)$ is a left $H^{\s}$-module map.
\end{lemma}

\begin{proof} {\rm (a)} It is clear that $1\rhu m=m$ for all $m\in M$. Now let $h$, $l\in H^{\s}$ and $m\in M$. Then we have
$$\begin{array}{rl}
& l\rhu(h\rhu m)\\
=&\S (l\rhu(h_{(2)}\cdot m_{(0)})_{(0)})
\s((h_{(2)}\cdot m_{(0)})_{(1)}\ot h_{(1)})\s^{-1}(h_{(3)}\ot m_{(1)})\\
=&\S(l_{(2)}\cdot(h_{(2)}\cdot m_{(0)})_{(0)})_{(0)}
\s((l_{(2)}\cdot(h_{(2)}\cdot m_{(0)})_{(0)})_{(1)}\ot l_{(1)})\\
&\s^{-1}(l_{(3)}\ot(h_{(2)}\cdot m_{(0)})_{(1)})
\s((h_{(2)}\cdot m_{(0)})_{(2)}\ot h_{(1)})\s^{-1}(h_{(3)}\ot m_{(1)})\\
\stackrel{(\ref{+-})}{=}&\S(l_{(2)}\cdot(h_{(3)}\cdot m_{(0)})_{(0)})_{(0)}
\s((l_{(2)}\cdot(h_{(3)}\cdot m_{(0)})_{(0)})_{(1)}\ot l_{(1)})\\
&\s(l_{(3)}(h_{(3)}\cdot m_{(0)})_{(1)}\ot h_{(1)})\s^{-1}(l_{(4)}\ot(h_{(3)}\cdot m_{(0)})_{(2)}h_{(2)})\\
&\s^{-1}(h_{(4)}\ot m_{(1)})\\
\stackrel{(2)}{=}&\S(l_{(2)}\cdot(h_{(2)}\cdot m_{(0)})_{(0)})_{(0)}
\s((l_{(2)}\cdot(h_{(2)}\cdot m_{(0)})_{(0)})_{(1)}\ot l_{(1)})\\
&\s(l_{(3)}(h_{(2)}\cdot m_{(0)})_{(1)}\ot h_{(1)})\s^{-1}(l_{(4)}\ot h_{(3)}m_{(1)})
\s^{-1}(h_{(4)}\ot m_{(2)})\\
\stackrel{(2)(\ref{rcocycle})}{=}&\S((l_{(3)}h_{(2)})\cdot m_{(0)})_{(0)}
\s(((l_{(3)}h_{(2)})\cdot m_{(0)})_{(1)}\ot l_{(1)})\\
&\s(((l_{(3)}h_{(2)})\cdot m_{(0)})_{(2)}l_{(2)}\ot h_{(1)})
\s^{-1}(l_{(4)}h_{(3)}\ot m_{(1)})\\
&\s^{-1}(l_{(5)}\ot h_{(4)})\\
\stackrel{(\ref{lcocycle})}{=}&\S((l_{(3)}h_{(3)})\cdot m_{(0)})_{(0)}
\s(l_{(1)}\ot h_{(1)})\s(((l_{(3)}h_{(3)})\cdot m_{(0)})_{(1)}\ot l_{(2)}h_{(2)})\\
&\s^{-1}(l_{(4)}h_{(4)}\ot m_{(1)})\s^{-1}(l_{(5)}\ot h_{(5)})\\
=&\S((l_{(2)}h_{(2)})\rhu m)\s(l_{(1)}\ot h_{(1)})\s^{-1}(l_{(3)}\ot h_{(3)})\\
=&(l\cdot_{\s}h)\rhu m.
\end{array}$$
This shows that $M$ is a left $H^{\s}$-module.

{\rm (b)} If $f:M\r N$ is a {\rm YD} $H$-module map, then $f$ is $H$-linear and $H$-colinear. By the definition of (\ref{sigmaaction}), we have for $h\in H^{\s}$ and $m\in M$,
\begin{eqnarray*}
f(h\rhu m)&=&\S f((h_{(2)}\cdot m_{(0)})_{(0)})\s((h_{(2)}\cdot m_{(0)})_{(1)}\ot h_{(1)}) \s^{-1}(h_{(3)}\ot m_{(1)})\\
&=&\S (h_{(2)}\cdot f(m)_{(0)})_{(0)})\s((h_{(2)}\cdot f(m)_{(0)})_{(1)}\ot h_{(1)}) \s^{-1}(h_{(3)}\ot f(m)_{(1)})\\
&=& h\rhu f(m).
\end{eqnarray*}
\end{proof}
Following \leref{2.1}, we obtain a covariant functor $\us$ from the category ${_H\mathcal{YD}}^H$  to the category of left $H^{\s}$-modules. We show that the functor $\us$ is in fact a covariant functor into the category  $\YDHs$.

\begin{lemma}\lelabel{2.2} Let $M$ be a {\rm YD} $H$-module. Then $\u{\s}(M)$ is a {\rm YD} $H^{\s}$-module with the inherited right $H$-comodule structure of $M$.
\end{lemma}

\begin{proof} Note that $H^{\s}=H$ as a coalgebra. Since $M$ is a right $H$-comodule, $\u{\s}(M)=M$ is also a right $H^{\s}$-comodule. By \leref{2.1}, $\u{\s}(M)$ is a left $H^{\s}$-module. It remains to show that $\us(M)$ satisfies the {\rm YD} $H^{\s}$-module compatibility (\ref{compat}).
Now let $h\in H^{\s}$ and $m\in\u{\s}(M)$. We have
$$\begin{array}{cl}
&\S(h_{(2)}\rhu m)_{(0)}\ot(h_{(2)}\rhu m)_{(1)}\cdot_{\s}h_{(1)}\\
=&\S(h_{(3)}\cdot m_{(0)})_{(0)}\ot(h_{(3)}\cdot m_{(0)})_{(1)}\cdot_{\s}h_{(1)}
\s((h_{(3)}\cdot m_{(0)})_{(2)}\ot h_{(2)})\s^{-1}(h_{(4)}\ot m_{(1)})\\
=&\S(h_{(3)}\cdot m_{(0)})_{(0)}\ot\s((h_{(3)}\cdot m_{(0)})_{(1)}\ot h_{(1)})
(h_{(3)}\cdot m_{(0)})_{(2)}h_{(2)}\s^{-1}(h_{(4)}\ot m_{(1)})\\
=&\S(h_{(2)}\cdot m_{(0)})_{(0)}\s((h_{(2)}\cdot m_{(0)})_{(1)}\ot h_{(1)})
\ot h_{(3)}m_{(1)}\s^{-1}(h_{(4)}\ot m_{(2)})\\
=&\S(h_{(2)}\cdot m_{(0)})_{(0)}\s((h_{(2)}\cdot m_{(0)})_{(1)}\ot h_{(1)})\s^{-1}(h_{(3)}\ot m_{(1)})\\
&\ot \s(h_{(4)}\ot m_{(2)})h_{(5)}m_{(3)}\s^{-1}(h_{(6)}\ot m_{(4)})\\
=&\S(h_{(1)}\rhu m_{(0)})\ot h_{(2)}\cdot_{\s}m_{(1)},
\end{array}$$
where the third equality follows from (\ref{compat}). It follows that
$\u{\s}(M)$ is a {\rm YD} $H^{\s}$-module.
\end{proof}

Thus we have constructed a covariant functor $\us$ from the category ${_{H}\mathcal{YD}}^{H}$ to the category $\YDHs$. We show that $\us$ is an equivalence monoidal functor and preserves the braidings. Consequently, the two braided monoidal categories ${_{H}\mathcal{YD}}^{H}$ and $\YDHs$ are  equivalent.

\begin{theorem}\thlabel{2.3} Let $\s$ be a $2$-cocycle on $H$.
Then $ \us~:{_H\mathcal{YD}}^H\r \YDHs$ is an equivalence braided monoidal functor.
\end{theorem}

\begin{proof} By \leref{2.1} and \leref{2.2}, $\us$ is a covariant functor. It remains to show that $\us$ is a braided monoidal functor and it has an inverse functor. For any {\rm YD} $H$-modules $M$ and $N$, we define an isomorphism
$\eta_{M,N}$:
$$ \u{\s}(M)\ot\u{\s}(N)\ra\u{\s}(M\ot N), m\ot n\mapsto
\S m_{(0)}\ot n_{(0)}\s^{-1}(n_{(1)}\ot m_{(1)}).$$
It is easy to see that $(\u{\s}, \eta)$ is a  monoidal functor from ${_H\mathcal{YD}}^H$ to $\YDHs$. We show that $(\u{\s}, \eta)$ commutes with the braidings in the two categories. Recall that the braiding $\Phi$ of $\YDH$ is defined by
$$\Phi_{M,N}: M\otimes N\r N\otimes M, \ \ \phi(m\otimes n)=\sum \no\ot \nl\cdot m$$
for $m\in M, n\in N$ and $M,N\in\; \YDH$. We have to verify  the following commutative diagram:
$$\begin{diagram}
\us(M)\otimes \us(N) & \rTo^{\Phi_{\us(M),\us(N)}} & \us(N)\otimes \us(M)\\
\dTo_{\eta_{M,N}} & & \dTo_{\eta_{N,M}}\\
\us(M\otimes N) & \rTo^{\us(\Phi_{M,N})} & \us(N\otimes M).
\end{diagram}$$
Let $m\in \us(M)$ and $n\in \us(N)$. Then we have
\begin{eqnarray*}
&& (\eta_{N,M}\Phi_{\us(M),\us(N)})(m\ot n) \\
&=& \sum \eta_{N,M}(\no\ot \nl\rhu m)\\
&=& \sum \eta_{N,M}(\no\ot (n_{(2)}\cdot m_{(0)})_{(0)})\s((n_{(2)}\cdot m_{(0)})_{(1)}\ot n_{(1)}) \s^{-1}(n_{(3)}\ot m_{(1)})\\
&=& \sum \no\ot (n_{(3)}\cdot m_{(0)})_{(0)} \s^{-1}((n_{(3)}\cdot m_{(0)})_{(1)}\ot \nl)
\s((n_{(3)}\cdot m_{(0)})_{(2)}\ot n_{(2)})\\
&& \s^{-1}(n_{(4)}\ot m_{(1)})\\
&=& \sum \no\ot (\nl\cdot m_{(0)})\s^{-1}(n_{(2)}\ot m_{(1)})\\
&=& (\us (\Phi_{M,N})\eta_{M,N})(m\ot n).
\end{eqnarray*}

Next we show that the braided monoidal functor $\us$ has an inverse functor. Recall from \cite[Lemma 1.2]{Ch0} that $\s^{-1}$ is a 2-cocycle on $H^{\s}$. It is easy to see that $(H^{\s})^{\s^{-1}}=H$ as Hopf algebras. Thus $\s^{-1}$ induces  a braided monoidal functor $(\us^{-1}, \xi)$ from $\YDHs$ to $\YDH$ with
$\us^{-1}(M)=M$ as right comodules, $\us^{-1}(f)=f$  and
$$\xi_{M,N}: \us^{-1}(M)\ot\us^{-1}(N)\ra \us^{-1}(M\ot N),\hspace{0.1cm}m\ot n\mapsto \S m_{(0)}\ot n_{(0)}\s(n_{(1)}\ot m_{(1)})$$
for all objects $M$, $N$ and all morphisms $f$ in the category $\YDHs$.

Now let $M$ be a {\rm YD} $H$-module. Denote by $\rightharpoondown$ the left
 $(H^{\s})^{\s^{-1}}=H$-action on $\us^{-1}(\u{\s}(M))$. Then for any $m\in\us^{-1}(\u{\s}(M))$
and $h\in H$, we have
$$\begin{array}{rcl}
h\rightharpoondown m&=&\S(h_{(2)}\rhu m_{(0)})_{(0)}\s^{-1}((h_{(2)}\rhu m_{(0)})_{(1)}\ot
h_{(1)})\s(h_{(3)}\ot m_{(1)})\\
&=&\S(h_{(3)}\cdot m_{(0)})_{(0)}\s((h_{(3)}\cdot m_{(0)})_{(2)}\ot h_{(2)})
\s^{-1}(h_{(4)}\ot m_{(1)})\\
&&\s^{-1}((h_{(3)}\cdot m_{(0)})_{(1)}\ot h_{(1)})
\s(h_{(5)}\ot m_{(2)})\\
&=&h\cdot m.
\end{array}$$
That is, $(\us^{-1}\hspace{0.1cm}\u{\s})(M)=M$. On the other hand,
it is not difficult to see that $\us^{-1}(\eta_{M,N})\circ\xi_{\u{\s}(M),\u{\s}(N)}$
is the identity on $M\ot N$, where we have $M\ot N=(\us^{-1}\hspace{0.1cm}\u{\s})(M)
\ot (\us^{-1}\hspace{0.1cm}\u{\s})(N)=(\us^{-1}\hspace{0.1cm}\u{\s})(M\ot N)$
for any {\rm YD} $H$-modules $M$ and $N$.
This shows that $(\us^{-1}, \xi)\circ(\u{\s}, \eta)$ is equal to the identity
functor $(\u{Id}, id)$ on the braided monoidal category $\YDH$.
Similarly, we have $(\u{\s}, \eta)\circ(\us^{-1}, \xi)=(\u{Id}, id)$ on the braided monoidal category $\YDHs$. Hence $\us$ is an equivalence braided monoidal functor. Consequently, $\YDHs$ and $\YDH$ are equivalent braided monoidal categories.
\end{proof}

Now we consider a CQT Hopf algebra $(H,R)$. Let $\s$ be a cocycle on $H$. By \cite{DT},  we know that $H^{\s}$ is also a CQT Hopf algebra with a CQT structure $R^{\s}=(\s\tau)*R*\s^{-1}$, that is,
$$R^{\s}(g\ot h)=\S\s(h_{(1)}\ot g_{(1)})R(g_{(2)}\ot h_{(2)})\s^{-1}(g_{(3)}\ot h_{(3)})$$
for all $g$, $h\in H$. Now let $M\in{\mathcal M}^H_R$. In this case, the $H$-module structure is given by (\ref{Raction}). We show that the YD $H^{\s}$-module $\us(M)$ is in ${\mathcal M}^{H^{\s}}_{R^{\s}}$. Namely, the $H^{\s}$-module structure of $\us(M)$ is induced by the CQT structure $R^{\s}$ through (\ref{Raction}). Thus we get the following well-known result (see the dual version \cite[Lemma XV 3.7]{Ka}).

\begin{corollary}\colabel{2.4}
Let $(H,R)$ be a CQT Hopf algebra and $\s$ be a cocycle on $H$. Then the equivalence braided monoidal functor $\us$ restricts to an equivalence braided monoidal  functor from the category $\M^H_R$ to $\M^{H^{\s}}_{R^{\s}}$.
\end{corollary}

\begin{proof}
Let $M\in \M^H_R$. To show $\us(M)\in \M^{H^{\s}}_{R^{\s}}$, we
have to verify that the $H^{\s}$-action (\ref{sigmaaction}) on
$\us(M)$ coincides with the $H^{\s}$-action, denoted
$\vartriangleright_1^{\s}$, induced by $R^{\s}$ and given by
(\ref{Raction}).  Indeed, let $m\in \us(M)$ and $h\in H^{\s}$. We
have
$$\begin{array}{rcl}
h\vartriangleright_1^{\s} m&=&\S m_{(0)}R^{\s}(h\ot m_{(1)})\\
&=&\S m_{(0)}\s(m_{(1)}\ot h_{(1)})R(h_{(2)}\ot m_{(2)})\s^{-1}(h_{(3)}\ot m_{(3)})\\
&=&\S m_{(0)}R(h_{(3)}\ot m_{(1)}) \s(h_{(4)}m_{(2)}S^{-1}(h_{(2)})\ot h_{(1)})
\s^{-1}(h_{(5)}\ot m_{(3)})\\
&=&\S(h_{(3)}\vartriangleright_1 m_{(0)})\s(h_{(4)}m_{(1)}S^{-1}(h_{(2)})\ot h_{(1)})
\s^{-1}(h_{(5)}\ot m_{(2)})\\
&=&h\rhu m,
\end{array}$$
 where we use (CQT$4''$) in the third equality. So  $\us(M)\in \M^{H^{\s}}_{R^{\s}}$.
\end{proof}

Let $A$ be a YD $H$-module algebra, namely, an algebra in the category $\YDH$.
Then $\u{\s}(A)$ is a {\rm YD} $H^{\s}$-module algebra with the multiplication given by
$$m^{\s}:\u{\s}(A)\ot\u{\s}(A)\stackrel{\eta_{A,A}}{\longrightarrow}
\u{\s}(A\ot A)\stackrel{\u{\s}(m)}{\longrightarrow}\u{\s}(A),$$
where $m:A\ot A\r A$ is the multiplication of $A$.
Denote by $a\bullet b$ the product $m^{\s}(a\ot b)$ for any $a$, $b\in\u{\s}(A)$. Then
\begin{equation}\label{sigmaproduct}
a\bullet b=\S a_{(0)}b_{(0)}\s^{-1}(b_{(1)}\ot a_{(1)}),\hspace{0.1cm}a, b\in\u{\s}(A).
\end{equation}
That is, $m^{\s}=m*(\s^{-1}\tau)$. Hence as an algebra $\u{\s}(A)$ is the same as $A_{\s^{-1}\tau}$ given in \cite[(2.27)]{Maj}, where $\tau$ is the usual flip map.

Recall from \cite{VZ0} that a braided monoidal functor $\mathcal{F}$ from a braided monoidal category $\mathcal{C}$ to a braided monoidal category $\mathcal{D}$ sends an Azumaya algebra in $\mathcal{C}$ to an Azumaya algebra in $\mathcal{D}$, and consequently induces a group homomorphism from the Brauer group Br$(\mathcal{C})$ to the Brauer group Br$(\mathcal{D})$. If the functor $\mathcal{F}$ is an equivalence functor, then the two Brauer groups are isomorphic. Thus the equivalence functor $\us$ from $\YDH$ to $\YDHs$ sends an $H$-Azumaya algebra $A$ to the $H^{\s}$-Azumaya algebra $\us(A)$ with the product (\ref{sigmaproduct}) and induces an group isomorphism from the Brauer group $\BQ(k,H)$ to the Brauer group $\BQ(k,H^{\s})$.

\begin{corollary}\colabel{2.5}
Let $H$ be a Hopf algebra and $\s$ be a  cocycle on $H$. Then the Brauer groups $\BQ(k,H)$ and $\BQ(k,H^{\s})$ are isomorphic. Furthermore, if $(H,R)$ is a CQT Hopf algebra, then the two equivariant Brauer groups $\BC(k,H,R)$ and $\BC(k,H^{\s},R^{\s})$ are isomorphic \cite{Car}.
\end{corollary}

In general, if $\s$ is a  2-cocycle on $H$ and $\s_1$ is a  2-cocycle on $H^{\s}$, then $\s_1*\s$ is a  2-cocycle on $H$ and $H^{\s_1*\s}=(H^{\s})^{\s_1}$ by \cite[Lemma 1.4]{Ch0}. Moreover, it is not difficult to check that
$\u{\s_1}(\u{\s}(A))=(\u{\s_1*\s})(A)$ as  {\rm YD}
$H^{\s_1*\s}$-modules (or module algebras) if $A$ is a {\rm YD} $H$-module (or module algebra).

A 2-cocycle $\s$ is called {\it lazy} if for all $h$, $l\in H$
$$\S\s(h_{(1)}\ot l_{(1)})h_{(2)}l_{(2)}=\S h_{(1)}l_{(1)}\s(h_{(2)}\ot l_{(2)}).$$
In other words, a lazy cocycle $\sigma$ commutes with the
multiplication of $H$. It is clear that a 2-cocycle $\s$ is lazy
if and only if $H^{\s}=H$ as Hopf algebras (\cite[Lemma
1.3]{Ch0}), i.e., the multiplication $\cdot_{\s}$ of $H^{\s}$ is
the same as the original multiplication of $H$. Thus, if $\s$ is a
lazy  2-cocycle on $H$ then $\s^{-1}$ is also a 2-cocycle on $H$
by \cite[Lemma 1.2]{Ch0}. The set of all lazy  2-cocycles forms a
group \cite{Ch0, BC}, denoted $\mathrm{Z}^2_L(H)$. The following
corollary tells us that the group $\mathrm{Z}^2_L(H)$ acts on the
Brauer group $\BQ(k,H)$ by automorphisms.

\begin{corollary}\colabel{2.6} Let $\s$ be a  $2$-cocycle on $H$. If $\s$ is lazy, then $\s$ induces an automorphism of the Brauer group ${\rm BQ}(k, H)$
$$\u{\s}: {\rm BQ}(k, H)\ra{\rm BQ}(k, H), [A]\mapsto [\u{\s}(A)].$$
\end{corollary}
\begin{proof} Follows \coref{2.5}.
\end{proof}

Note that in \cite{Sch} P. Schauenburg introduced a cohomology group $\mathrm{H}^2_c(H)$  which is the quotient group of $\mathrm{Z}^2_L(H)$ modulo the  subgroup $\mathrm{B}^2_L(H)$ consisting of  coboundary lazy cocycles. The cohomology group has been systematically studied by J. Bichion and G. Carnovale in \cite{BC}, where the group is called the lazy cohomology group, denoted $\mathrm{H}^2_L(H)$. Precisely, a  lazy cocycle $\sigma$ is coboundary if there is an (invertible) lazy 1-cocycle $\mu:H\r k$ such that $\sigma(a,b)=\mu(a_{(1)})\mu(\bl)\mu^{-1}(\az\bz)$ for all $a,b\in H$, where $\mu$ satisfies the relations:
 $\sum \mu(\al)\az=\sum \mu(\az)\al$ and $\mu(1)=1$ for all $a\in H$. That is, $\mu$ is a normalized central element in $H^*$. We show that if $\sigma$ is a  coboundary lazy cocycle, then the functor $\us$ is isomorphic to the identity functor of $\YDH$. Thus the subgroup $\mathrm{B}^2_L(H)$ acts on $\BQ(k,H)$ trivially. Consequently, the lazy cohomology group $\mathrm{H}^2_L(H)$ defined in \cite{Sch,BC} acts on $\BQ(k,H)$ by automorphisms.

\begin{corollary}\colabel{2.7} Let $\s$ be a coboundary lazy $2$-cocycle on $H$, then $\us$ is isomorphic to the identity functor. Consequently, the lazy cohomology group $\mathrm{H}^2_L(H)$  acts on $\BQ(k,H)$ by automorphisms.
\end{corollary}

\begin{proof}
Let $\sigma\in \mathrm{B}^2_L(H)$. That is,  $\sigma(g\otimes h)=\sum\mu(g_{(1)})\mu(h_{(1)})\mu^{-1}(g_{(2)}h_{(2)})$ for any $g, h\in H$, where $\mu\in H^*$ is an invertible element such that
$\mu(1)=1$ and $\sum\mu(h_{(1)})h_{(2)}=\sum h_{(1)}\mu(h_{(2)})$
for any $h\in H$. It is clear that $\sigma^{-1}(g\otimes
h)=\sum\mu(g_{(1)}h_{(1)})\mu^{-1}(g_{(2)})\mu^{-1}(h_{(2)})$ for
any $g, h\in H$.

Now for any $M\in{\ }_H{\mathcal YD}^H$, define
$$\zeta_M: M\rightarrow\underline{\sigma}(M),\ \ m\mapsto\sum
m_{(0)}\mu(m_{(1)}).$$
Then $\zeta_M$ is a YD $H$-module isomorphism. In fact, since $\mu$ is invertible, $\zeta_M$ is a $k$-linear isomorphism. Now for any $m\in M$ and $h\in H$, we have
$$\begin{array}{rcl}
\rho(\zeta_M(m))&=&\sum\rho(m_{(0)})\mu(m_{(1)})\\
&=&\sum m_{(0)}\otimes m_{(1)}\mu(m_{(2)})\\
&=&\sum m_{(0)}\otimes\mu(m_{(1)}) m_{(2)}\\
&=&\sum\zeta_M(m_{(0)})\otimes m_{(1)},\\
\end{array}$$
and
$$\begin{array}{rcl}
h\rightharpoonup\zeta_M(m)&=&\sum(h\rightharpoonup
m_{(0)})\mu(m_{(1)})\\
&=&\sum(h_{(2)}\cdot m_{(0)})_{(0)}\sigma((h_{(2)}\cdot
m_{(0)})_{(1)}\otimes h_{(1)})\sigma^{-1}(h_{(3)}\otimes
m_{(1)})\mu(m_{(2)})\\
&=&\sum(h_{(3)}\cdot m_{(0)})_{(0)}\mu((h_{(3)}\cdot
m_{(0)})_{(1)})\mu(h_{(1)})\mu^{-1}((h_{(3)}\cdot
m_{(0)})_{(2)}h_{(2)})\\
& &\mu(h_{(4)}m_{(1)})\mu^{-1}(h_{(5)})\mu^{-1}(m_{(2)})\mu(m_{(3)})\\
&=&\sum(h_{(3)}\cdot m_{(0)})_{(0)}\mu((h_{(3)}\cdot
m_{(0)})_{(1)})\mu(h_{(1)})\mu^{-1}((h_{(3)}\cdot
m_{(0)})_{(2)}h_{(2)})\\
& &\mu(h_{(4)}m_{(1)})\mu^{-1}(h_{(5)})\\
&=&\sum(h_{(2)}\cdot m_{(0)})_{(0)}\mu((h_{(2)}\cdot
m_{(0)})_{(1)})\mu^{-1}((h_{(2)}\cdot m_{(0)})_{(2)}h_{(1)})\\
& &\mu(h_{(3)}m_{(1)})\ \ \ (\mbox{ since }\sum\mu(h_{(1)})h_{(2)}=\sum h_{(1)}\mu(h_{(2)}))\\
&=&\sum\zeta_M((h_{(2)}\cdot m_{(0)})_{(0)})\mu^{-1}((h_{(2)}\cdot
m_{(0)})_{(1)}h_{(1)})\mu(h_{(3)}m_{(1)})\\
&\stackrel{(2)}{=}&\sum\zeta_M(h_{(1)}\cdot
m_{(0)})\mu^{-1}(h_{(2)}m_{(1)})\mu(h_{(3)}m_{(2)})\\
&=&\zeta_M(h\cdot m).\\
\end{array}$$
This shows that $\zeta_M$ is a YD $H$-module isomorphism from $M$
to $\underline{\sigma}(M)$. It is easy to see that $\zeta$ is a
natural transformation from the identity functor $Id$ to the functor
$\underline{\sigma}$. Furthermore, for any YD $H$-modules $M$ and $N$, we have the
following commutative diagram:
$$\begin{array}{rcl}
&M\otimes N&\\
\zeta_M\otimes\zeta_N\swarrow& &\searrow\zeta_{M\otimes N}\\
\underline{\sigma}(M)\otimes\underline{\sigma}(N)&
\stackrel{\eta_{M, N}}{\longrightarrow}&
\underline{\sigma}(M\otimes N).\\
\end{array}$$
In fact, let $m\in M$ and $n\in N$. We have
$$\begin{array}{rcl}
(\eta_{M,N}(\zeta_M\otimes\zeta_N))(m\otimes n)&=&
\sum\eta_{M,N}(m_{(0)}\otimes n_{(0)})\mu(m_{(1)})\mu(n_{(1)})\\
&=&\sum m_{(0)}\otimes n_{(0)}\sigma^{-1}(n_{(1)}\otimes m_{(1)})\mu(m_{(2)})\mu(n_{(2)})\\
&=&\sum m_{(0)}\otimes n_{(0)}\mu(n_{(1)}m_{(1)})\\
&=&\sum (m\otimes n)_{(0)}\mu((m\otimes n)_{(1)})\\
&=&\zeta_{M\otimes N}(m\otimes n).\\
\end{array}$$
It follows from the above commutative diagram that $\zeta_A:
A\rightarrow\underline{\sigma}(A)$ is also an algebra
homomorphism, and hence a YD $H$-module algebra isomorphism if $A$
is a YD $H$-module algebra. Thus we have proved that $\mathrm{B}^2_L(H)$ acts trivially
on the Brauer group BQ$(k, H)$. Hence $\mathrm{H}^2_L(H)$ acts on BQ$(k,
H)$ by automorphisms.
\end{proof}

For completeness, we now give the dual versions of Theorem 2.3 and Corollary 2.4. We will omit the detail of proofs and give the sketch of the construction.
Let $H$ be a Hopf algebra with a bijective antipode in the sequel.
For an element $r\in H\ot H$, we let $r_{12}$ and $r_{23}$ stand for $r\ot 1$ and $ 1\ot r$ in $H\ot H\ot H$ respectively.

Recall from \cite{Maj} that an invertible element $\theta=\sum \theta^{(1)}\ot \theta^{(2)}\in H\ot H$ is called a {\it dual $2$-cocycle} if it satisfies:
$$\theta_{12}(\D\ot{\rm id})(\theta)=\theta_{23}({\rm id}\ot\D)(\theta).$$
and  $(\e\ot{\rm id})({\theta})=({\rm id}\ot\e)(\theta)=1$.

To a dual  2-cocycle $\theta\in H\ot H$ one may associate a new Hopf algebra $H_{\theta}$. As an algebra $H_{\theta}=H$ but with comultiplication given by
$$\D_{\theta}(h)=\theta\D(h)\theta^{-1},\hspace{0.2cm}h\in\hspace{0.1cm}H_{\theta}.$$
The antipode $S_{\theta}$ of $H_{\theta}$ is given by
\begin{eqnarray*}
S_{\theta}(h)&=&\S\theta^{(1)}S(\theta^{(2)})S(h)S((\theta^{-1})^{(1)})
(\theta^{-1})^{(2)}\\
&=&\S \theta^{(1)}S((\theta^{-1})^{(1)}h\theta^{(2)})(\theta^{-1})^{(2)},
\hspace{0.2cm}h\in\hspace{0.1cm}H_{\theta},
\end{eqnarray*}
where $\theta=\S \theta^{(1)}\ot\theta^{(2)}$ and $\theta^{-1}=\S(\theta^{-1})^{(1)}\ot(\theta^{-1})^{(2)}$
in $H\ot H$. Note that $S_{\theta}$ is bijective since $S$ is bijective.

Let $M$ be a Yetter-Drinfeld $H$-module with the coaction $\rho(m)=\S m_{(0)}\ot m_{(1)}$, $m\in M$.
Then one can define a new $H_{\theta}$-coaction
$\rho_{\theta}: M\ra M\ot\hspace{0.1cm}H_{\theta}$
by
\begin{equation}\label{newcom}
\begin{array}{rcl}
\hspace{-0.2cm}\rho_{\theta}(m)&=&\S\theta^{(1)}\cdot((\theta^{-1})^{(2)}\cdot m)_{(0)}
\ot\theta^{(2)}((\theta^{-1})^{(2)}\cdot m)_{(1)}(\theta^{-1})^{(1)}\\
&=&\S(\theta^{(1)}(\theta^{-1})^{(2)}_{(2)})\cdot m_{(0)}\ot
\theta^{(2)}(\theta^{-1})^{(2)}_{(3)}m_{(1)}S^{-1}((\theta^{-1})^{(2)}_{(1)})(\theta^{-1})^{(1)}.
\end{array}
\end{equation}
Since $H_{\theta}=H$ as algebras, $M$ is a left $H_{\theta}$-module.
It is straightforward to verify that the inherited $H_{\theta}$-module $M$ with  the $H_{\theta}$-comodule structure $\rho_{\theta}$ given by (\ref{newcom}) is a Yetter-Drinfeld $H_{\theta}$-module, denoted $\underline{\theta}(M)$. Thus a  dual cocycle $\theta$ induces a covariant functor $\underline{\theta}$ from $_H\mathcal{YD}^H$ to $_{H_{\theta}}\mathcal{YD}^{H_{\theta}}$ sending a YD $H$-module
 $M$ to the YD $H_{\theta}$-module $\underline{\theta}(M)$ and a YD $H$-module morphism  $f$ to a YD $H_{\theta}$-module morphism $\underline{\theta}(f)=f$.

For any two Yetter-Drinfeld $H$-modules $M$ and $N$, let
$$\phi_{M, N}: \underline{\theta}(M)\ot\underline{\theta}(N)\ra\underline{\theta}(M\ot N),
\hspace{0.2cm}m\ot n\mapsto\theta^{-1}\cdot(m\ot n)=\S(\theta^{-1})^{(1)}\cdot m\ot(\theta^{-1})^{(2)}\cdot n.$$
Then $\phi_{M,N}$ is a family of natural isomorphisms and
$(\underline{\theta}, \phi)$ gives a braided monoidal category equivalence
from $_H\mathcal{YD}^H$ to $_{H_{\theta}}\mathcal{YD}^{H_{\theta}}$. We summarize our above discussion in the following dual version of \thref{2.3}.

\begin{theorem}\thlabel{2.8}
Let $\theta\in H\ot H$ be a dual $2$-cocycle. Then the functor $(\underline{\theta},\phi)$ is an equivalence braided monoidal functor from $_H\mathcal{YD}^H$ to $_{H_{\theta}}\mathcal{YD}^{H_{\theta}}$.
\end{theorem}

Now assume that $H$ is a QT Hopf algebra with a QT structure
$\R=\S \R^{(1)}\ot \R^{(2)}$. Let
$$\R_{\theta}:=\tau(\theta)\R\theta^{-1}=\S \theta^{(2)}\R^{(1)}(\theta^{-1})^{(1)}\ot\theta^{(1)}\R^{(2)}(\theta^{-1})^{(2)}.$$
Then $\R_{\theta}$ is a QT structure of $H_{\theta}$.
In this case, restricting the functor $(\underline{\theta}, \phi)$ on the subcategory
$_H\!\mathcal{M}^\R$ gives a braided monoidal category equivalence from $_H\!\mathcal{M}^\R$
to $_{H_{\theta}}\!\mathcal{M}^{\R_{\theta}}$. Hence we have the following corollary dual to \coref{2.4}.

\begin{corollary}\colabel{2.9} \cite[Lemma XV 3.7]{Ka}
Let $(H, \R)$ be a QT Hopf algebra and let $\theta\in H\ot H$ be a  dual  $2$-cocycle. Then the equivalence braided monoidal functor $(\u{\theta},\phi)$ restricts to an equivalence braided monoidal functor from $_{H}\M^\R$ to $_{H_{\theta}}\M^{\R_{\theta}}$.
\end{corollary}

Let $A$ be a Yetter-Drinfeld $H$-module algebra with multiplication $m: A\ot A\ra A$. It follows from \thref{2.8} that $\underline{\theta}(A)$ is a Yetter-Drinfeld  $H_{\theta}$-module algebra with the multiplication given by
$$m_{\theta}: \underline{\theta}(A)\ot\underline{\theta}(A)
\stackrel{\phi_{A,A}}{\longrightarrow}\underline{\theta}(A\ot A)
\stackrel{\underline{\theta}(m)}{\longrightarrow}\underline{\theta}(A).$$
Let $a\bullet b:=m_{\theta}(a\ot b)$ for any $a$, $b\in\underline{\theta}(A)$. Then
$$a\bullet b=\S((\theta^{-1})^{(1)}\cdot a)((\theta^{-1})^{(2)}\cdot b),
\hspace{0.2cm}a, b\in\underline{\theta}(A).$$
If $A$ is an $H$-Azumaya algebra, then $\u{\theta}(A)$ is an $H_{\theta}$-Azumaya algebra. Thus, the functor $\u{\theta}$ induces an group isomorphism from BQ$(k,H)$ to BQ$(k,H_{\theta})$.

\begin{corollary}\colabel{2.10}
Let $H$ be a Hopf algebra and $\theta$ be a dual $2$-cocycle in $H$.
Then the Brauer group $\BQ(k,H)$ is isomorphic to the Brauer group $\BQ(k,H_{\theta})$. Moreover, if $(H,\R)$ is a QT Hopf algebra, then the two equivariant Brauer groups $\BM(k,H,\R)$ and $\BM(k,H_{\theta},\R_{\theta})$ are isomorphic.
\end{corollary}

A dual  2-cocycle $\theta\in H\ot H$ is said to be {\it lazy} if $\theta\D(h)=\D(h)\theta$ for all
$h\in H$. This is equivalent to $\D_{\theta}=\D$, i.e., $H_{\theta}=H$ as Hopf algebras. Note that for a lazy dual 2-cocycle $\theta$, its inverse $\theta^{-1}$ is a dual 2-cocycle as well. Let $\mathrm{Z}_2^L(H)$ be the set of all lazy dual 2-cocycles in $H$. Then $\mathrm{Z}_2^L(H)$
forms a group, which acts on the Brauer group $\BQ(k,H)$ by automorphisms like the group $\mathrm{Z}^2_L(H)$ does.  Let $\mathrm{B}^L_2(H)$ be the subgroup of $\mathrm{Z}^L_2(H)$ consisting of coboundary lazy dual cocycle $\theta$, i.e., $\theta=(u\otimes u)\Delta(u^{-1})$ for some central invertible element $u$ in $H$ with $\varepsilon(u)=1$. Similar to \coref{2.7}, the subgroup $\mathrm{B}^L_2(H)$ acts on $\BQ(k,H)$ trivially.

\begin{corollary}\colabel{2.11}
The cohomology group $\mathrm{H}^L_2(H)$ acts on the Brauer group $\BQ(k,H)$ by automorphisms.
\end{corollary}

\begin{example}\exlabel{2.12}
Let $k$ be a field with $ch(k)\not=2$. Let $H_4$ be Sweedler's 4-dimensional Hopf algebra over $k$ generated by two elements $g$ and $h$ satisfying the relations: $g^2=1,\ h^2=0,\ gh+hg=0$. The comultiplication, the counit and the antipode of $H_4$ are given by
$$\Delta(g)=g\ot g,\ \Delta(h)=1\ot h+h\ot g,\ \varepsilon(g)=1,\ \varepsilon(h)=0, \ S(g)=g, \ S(h)=gh.$$
It is well-known that $H_4$ has a family of QT structures parameterized by $t\in k$:
$$\R_t=\frac{1}{2}(1\ot 1+1\ot g+g\ot 1-g\ot g)+\frac{t}{2}(1\ot 1+g\ot g+1\ot
g-g\ot 1)(h\ot h).$$
A straightforward computation shows that the group $\mathrm{Z}_2^L(H_4)$ of the  lazy dual 2-cocycles $\theta_t$ is isomorphic to the additive group $k^+$, where $\theta_t$ are given as follows:
$$\theta_t=1\ot 1+\frac{t}{2}h\ot gh, \ \ t\in k.$$
Note that $H_4$ is a self-dual Hopf algebra. It has a family of CQT structures $R_t$ and the cohomology group $\mathrm{H}^2_L(H_4)=\mathrm{Z}^2_L(H_4)$ is isomorphic to $k^+$ \cite[Example 2.1]{BC}, where $R_t$ and lazy cocycles $\s_t$ are given as follows:
$$\begin{tabular}{l|llll}
$R_t$ & 1 & $g$ & $h$& $gh$\\ \hline
1 & 1 & 1 & 0 & 0\\
$g$ & 1& -1 & 0 & 0\\
$h$& 0& 0& $t$ & -$t$\\
$gh$& 0 & 0& $t$ & $t$,
\end{tabular}\  \ \ \ \ \hspace{1cm}
\begin{tabular}{l|llll}
$\sigma_t$ & 1 & $g$ & $h$& $gh$\\ \hline
1 & 1 & 1 & 0 & 0\\
$g$ & 1& 1 & 0 & 0\\
$h$& 0& 0& $t/2$ & -$t/2$\\
$gh$& 0 & 0& $t/2$ & -$t/2$,
\end{tabular}$$
where $t\in k$. The action of the automorphism $\u{\s_s}$ on $\BQ(k,H_4)$ induced by a lazy cocycle $\s_s$  is of particular interest, for $\u{\s_s}$ moves the equivariant Brauer group $\BC(k,H_4,R_t)$ to $\BC(k,H_4,R_{t-s})$ and fixes the subgroup $\BW(k)$, the Brauer-Wall group of $k$. By \cite{VZ1, Z} and \coref{2.5}, we have $\BC(k,H_4,R_t)\cong\BW(k)\times k_t^+$, where $k_t^+\cong k^+$. It is not hard to show that $k_t^+\cap k_s^+=1$ if $s\not=t$ in $k$. Thus, no distinct two elements $\u{\s_s}$ and $\u{\s_{s'}}$ act as the same automorphism of $\BQ(k,H_4)$. It follows that $\mathrm{H}^2_L(H_4)\cong k^+$ is a subgroup of the automorphism group of $\BQ(k,H_4)$.

Similarly, the group $\mathrm{H}_2^L(H_4)\cong k^+$ is a subgroup of the automorphism group of $\BQ(k,H_4)$. We show that the intersection of the two subgroups $\mathrm{H}^2_L(H_4)$ and $\mathrm{H}_2^L(H_4)$ is trivial. One may take a while to compute that $\u{\s_s}$ restricts to the identity automorphism of the subgroup $\BM(k,H_4,\R_0)$ for any $s\in k$. But $\u{\theta_{s'}}$ moves $\BM(k,H_4,\R_0)$ to the subgroup $\BM(k,H_4,\R_{-s'})$ for $s'\in k$. Thus $\u{\theta_s}\not=\u{\s_{s'}}$ as automorphisms except for $s=s'=0$.
\end{example}

\section{\bf The group of bigalois objects}\selabel{3}
\def\GalHR{\Gal(\HR)}
\def\GalHRs{\Gal({\mathcal H}^{\s}_{R^{\s}})}
Throughout this section $(H, R)$ will be a finite CQT Hopf algebra. In \cite{Z}, the author constructed a group $\Gal(\HR)$ of $\HR$-bigalois objects, where $\HR$ is a braided Hopf algebra in the category $\M^H_R$, in order to compute the equivariant Brauer group ${\rm BC}(k, H, R)$. In this section, we study the group $\Gal(\HR)$ under a cocycle deformation. We will show that the group ${\rm Gal}({\mathcal H}_R)$  is stable under a cocycle deformation. That is, $\Gal(\HR)$ and $\GalHRs$ are isomorphic groups.

Recall that ${\mathcal H}_R$ is a braided Hopf algebra in the category ${\mathcal M}^H_R$ defined by means of $H$. As a coalgebra, ${\mathcal H}_R$ coincides with $H$. The multiplication $\star$ and the antipode $S_R$ are given by
$$h\star l=\S l_{(2)}h_{(2)}R(S^{-1}(l_{(3)})l_{(1)}\ot h_{(1)}), \ {\rm and}$$
$$S_R(h)=\S S(h_{(2)})R(S^2(h_{(3)})S(h_{(1)})\ot h_{(4)})$$
respectively, where $h, l\in H$. As an object in ${\mathcal M}^H_R$, ${\mathcal H}_R$ has the adjoint coaction:
$$\rho(h)=\S h_{(2)}\ot S(h_{(1)})h_{(3)},\hspace{0.2cm}h\in{\mathcal H}_R.$$
For the detail of $\HR$, the reader can refer to \cite[Theorem 4.1]{Maj1} or \cite[Lemma 2.1]{Z}.

Let $M$ be a {\rm YD} $H$-module. We use the Sweedler notation $\S m_{(0)}\ot m_{(1)}$ for the $H$-comodule structure of $M$ as usual, and use the following summation notation for the dual $H^*$-comodule structure of the  $H$-module structure of $M$:
\begin{equation}\label{rcom}
M\ra M\ot H^*,\hspace{0.2cm}m\mapsto\S m_{[0]}\ot m_{[1]}.
\end{equation}
A natural left ${\mathcal H}_R$-module structure stemming from the YD $H$-module structure of $M$ is as follows:
\begin{equation}\label{LM}
\begin{array}{rcl}
h-\!\!\triangleright m&=&\S S^{-1}(h_{(2)})\vartriangleright_1(h_{(1)}\cdot m)\\
&=&\S(h_{(2)}\cdot m_{(0)})R(S^{-1}(h_{(4)})\ot h_{(3)}m_{(1)}S^{-1}(h_{(1)}))
\end{array}
\end{equation}
for $h\in{\mathcal H}_R$, $m\in M$, where $\vartriangleright_1$ stands for the action  (\ref{Raction}). Observe that the right $H$-comodule structure of $M$ induces two left $H$-module structures.
The first one is given by (\ref{Raction}), and the second one is given by
\begin{equation}\label{Lmod2}
h\vartriangleright_2 m=\S m_{(0)}R(S(m_{(1)})\ot h)=\S m_{(0)}R(m_{(1)}\ot S^{-1}(h))
\end{equation}
for $h\in H$ and $m\in M$. Note that a right $H$-comodule $M$ is simultaneously a
{\rm YD} $H$-module with the left $H$-action (\ref{Lmod2}). Now we define a right ${\mathcal H}_R$-module
structure on $M$ as follows:
\begin{equation}\label{RM}
\begin{array}{rcl}
m\triangleleft\!\!-h&=&\S S(h_{(1)})\vartriangleright_2(h_{(2)}\cdot m)\\
&=&\S(h_{(3)}\cdot m_{(0)})R(h_{(4)}m_{(1)}S^{-1}(h_{(2)})\ot h_{(1)})
\end{array}
\end{equation}
for $m\in M$ and $h\in{\mathcal H}_R$. It follows from \cite{Z} that $M$ is an
${\mathcal H}_R$-bimodule with the actions (\ref{LM}) and (\ref{RM}), and consequently
$M$ is an ${\mathcal H}_R^*$-bicomodule. Denote the dual left and right
${\mathcal H}_R^*$-comodule structures by
\begin{eqnarray*}
M\ra {\mathcal H}_R^*\ot M,&& m\mapsto\S m^{(-1)}\ot m^{(0)}, \ {\mathrm{and}}\\
M\ra M\ot{\mathcal H}_R^*,&& m\mapsto\S m^{(0)}\ot m^{(1)},
\end{eqnarray*}
respectively. Then $m\triangleleft\!\!-h=\S\langle m^{(-1)},h\rangle m^{(0)}$ and
$h-\!\!\triangleright m=\S\langle m^{(1)},h\rangle m^{(0)}$.

Similarly, we have a braided Hopf algebra ${\mathcal H}^{\s}_{R^{\s}}$
in the category ${\mathcal M}^{H^{\s}}_{R^{\s}}$, and one can define an
${\mathcal H}^{\s}_{R^{\s}}$-bimodule structure on any YD $H^{\s}$-module $M$. Consequently, $M$ has an
$({\mathcal H}^{\s}_{R^{\s}})$-bimodule structure similar to (\ref{LM}) and
(\ref{RM}) respectively. Denote the left and right dual $({\mathcal H}^{\s}_{R^{\s}})^*$-comodule structures by
\begin{eqnarray*}
M\ra ({\mathcal H}^{\s}_{R^{\s}})^*\ot M,&& m\mapsto\S m^{<-1>}\ot m^{<0>}, \ \mathrm{and}\\
M\ra M\ot({\mathcal H}^{\s}_{R^{\s}})^*,&& m\mapsto\S m^{<0>}\ot m^{<1>},
\end{eqnarray*}
respectively.  If $M$ is a {\rm YD} $H^{\s}$-module, then the left $H^{\s}$-module structure of $M$ induces a dual right $(H^{\s})^*$-comodule structure similar to (\ref{rcom}).
Let us denote the right $(H^{\s})^*$-comodule structure by
$$
M\ra M\ot(H^{\s})^*,\hspace{0.2cm}m\mapsto\S m_{<0>}\ot m_{<1>}.
$$

Now let $_{\diamond}M$ and $M_{\diamond}$ stand for the left and right
${\mathcal H}_R^*$-coinvariant submodules respectively. That is,
$$\begin{array}{rcl}
_{\diamond}M&=&\left\{m\in M| \S m^{(-1)}\ot m^{(0)}=\e\ot m\right\}\\
&=&\left\{m\in M| m\triangleleft\!\!-h=\e(h)m, \forall h\in{\mathcal H}_R\right\}
\end{array}$$ and
$$\begin{array}{rcl}
M_{\diamond}&=&\left\{m\in M| \S m^{(0)}\ot m^{(1)}=m\ot\e\right\}\\
&=&\left\{m\in M| h-\!\!\triangleright m=\e(h)m, \forall h\in{\mathcal H}_R\right\}.
\end{array}$$

\begin{lemma}\lelabel{3.1}\cite[Lemma 2.5, p332]{Z}
Let $M$ be a {\rm YD} $H$-module. Then we have\\
{\rm (a)}\hspace{0.3cm}$M_{\diamond}=\left\{ m\in M| h\cdot m=h\vartriangleright_1m, \forall h\in H\right\},$\\
{\rm (b)}\hspace{0.3cm}$_{\diamond}M=\left\{m\in M|h\cdot m=h\vartriangleright_2m, \forall h\in H\right\}.$
\end{lemma}
The following lemma says that the left and right $\HR$-invariant functors $(-)_{\diamond}$ and $_{\diamond}(-)$ are subfunctors from $\YDH$ to $\YDH$.

\begin{lemma}\lelabel{3.2}
Let $M$ be a {\rm YD} $H$-module. Then $M_{\diamond}$ and $_{\diamond}M$ are {\rm YD} $H$-submodule of $M$.
\end{lemma}

\begin{proof} Let $m\in M_{\diamond}$. It follows from \leref{3.1}(a) that
$h\cdot m=h\vartriangleright_1m=\S m_{(0)}R(h\ot m_{(1)})$
for any $h\in H$. Hence we have
$$\begin{array}{rcl}
\S h\cdot m_{(0)}\ot m_{(1)}&=&\S h_{(1)}\cdot m_{(0)}\ot S^{-1}(h_{(3)})h_{(2)}m_{(1)}\\
&=&\S(h_{(2)}\cdot m)_{(0)}\ot S^{-1}(h_{(3)})(h_{(2)}\cdot m)_{(1)}h_{(1)}\\
&=&\S (h_{(2)}\vartriangleright_1m)_{(0)}\ot S^{-1}(h_{(3)})(h_{(2)}\vartriangleright_1m)_{(1)}h_{(1)}\\
&=&\S m_{(0)}\ot S^{-1}(h_{(3)})m_{(1)}h_{(1)}R(h_{(2)}\ot m_{(2)})\\
&=&\S m_{(0)}\ot R(h_{(1)}\ot m_{(1)})S^{-1}(h_{(3)})h_{(2)}m_{(2)}\\
&=&\S m_{(0)}R(h\ot m_{(1)})\ot m_{(2)}\\
&=&\S h\vartriangleright_1m_{(0)}\ot m_{(1)}
\end{array}$$
for all $h\in H$. It follows from \leref{3.1}(a) that
$\S m_{(0)}\ot m_{(1)}\in M_{\diamond}\ot H$. This shows that $M_{\diamond}$ is an $H$-subcomodule of $M$.
Then for any $m\in M_{\diamond}$ and $h\in H$, we have
$h\cdot m=h\vartriangleright_1m=\S m_{(0)}R(h\ot m_{(1)})\in M_{\diamond}$.
Hence $M_{\diamond}$ is an $H$-submodule of $M$, and so $M_{\diamond}$ is a {\rm YD} $H$-submodule of $M$.
Similarly, using \leref{3.1}(b), one can show that $_{\diamond}M$ is also a {\rm YD} $H$-submodule of $M$.
\end{proof}
Next we show that the equivalence functor $\us$ commutes with both $\HR$-invariant functors $_{\diamond}(-)$ and $(-)_{\diamond}$.

\begin{lemma}\lelabel{3.3}
Let $M$ be a {\rm YD} $H$-module. Then $\u{\s}(M)_{\diamond}=M_{\diamond}$ and
$_{\diamond}\u{\s}(M)=\hspace{0.05cm}_{\diamond}M$ as $k$-modules, and hence
$\u{\s}(M_{\diamond})=\u{\s}(M)_{\diamond}$ and
$\u{\s}(_{\diamond}M)=\hspace{0.05cm}_{\diamond}\u{\s}(M)$ as {\rm YD} $H^{\s}$-modules. In particular, if $A$ is a {\rm YD} $H$-module algebra, then
$\u{\s}(A_{\diamond})=\u{\s}(A)_{\diamond}$ and
$\u{\s}(_{\diamond}A)=\hspace{0.05cm}_{\diamond}\u{\s}(A)$ as {\rm YD} $H^{\s}$-module algebras.
\end{lemma}

\begin{proof} By \coref{2.4} and its proof, one can see
$M_{\diamond}\subseteq\u{\s}(M)_{\diamond}$. Since $\u{\s}(M)$ is a {\rm YD} $H^{\s}$-module
and $\s^{-1}$ is a  2-cocycle on $H^{\s}$, we also have
$\u{\s}(M)_{\diamond}\subseteq\us^{-1}(\u{\s}(M))_{\diamond}=M_{\diamond}$.
Hence $\u{\s}(M)_{\diamond}=M_{\diamond}$.

Let $m\in\hspace{0.05cm}_{\diamond}M$ and $h\in H^{\s}$. Then it follows from
\leref{3.1} and \leref{3.2} that
$$\begin{array}{rcl}
h\vartriangleright_2^{\s}m&=&\S m_{(0)}R^{\s}(m_{(1)}\ot(S^{\s})^{-1}(h))\\
&=&\S m_{(0)}R^{\s}(m_{(1)}\ot S^{-1}(h_{(3)}))\s(S^{-1}(h_{(2)})\ot h_{(1)})\\
&&\s^{-1}(h_{(5)}\ot S^{-1}(h_{(4)}))\\
&=&\S m_{(0)}\s(S^{-1}(h_{(5)})\ot m_{(1)})R(m_{(2)}\ot S^{-1}(h_{(4)}))\\
&&\s^{-1}(m_{(3)}\ot S^{-1}(h_{(3)}))\s(S^{-1}(h_{(2)})\ot h_{(1)})\s^{-1}(h_{(7)}\ot S^{-1}(h_{(6)}))\\
&=&\S m_{(0)}R(m_{(1)}\ot S^{-1}(h_{(6)}))
\s(S^{-1}(h_{(8)})\ot h_{(7)}m_{(2)}S^{-1}(h_{(5)}))\\
&&\s(m_{(3)}S^{-1}(h_{(4)})\ot h_{(1)})\s^{-1}(m_{(4)}\ot S^{-1}(h_{(3)})h_{(2)})\\
&&\s^{-1}(h_{(10)}\ot S^{-1}(h_{(9)}))\\
&=&\S(h_{(4)}\vartriangleright_2m_{(0)})\s(S^{-1}(h_{(8)})\ot h_{(5)}m_{(1)}S^{-1}(h_{(3)}))\\
&&\s(S^{-1}(h_{(7)})h_{(6)}m_{(2)}S^{-1}(h_{(2)})\ot h_{(1)})\s^{-1}(h_{(10)}\ot S^{-1}(h_{(9)}))\\
&=&\S(h_{(3)}\cdot m_{(0)})\s(h_{(4)}m_{(1)}S^{-1}(h_{(2)})\ot h_{(1)})\\
&&\s(S^{-1}(h_{(6)})\ot h_{(5)}m_{(2)})\s^{-1}(h_{(8)}\ot S^{-1}(h_{(7)}))\\
&=&\S(h_{(3)}\cdot m_{(0)})\s(h_{(4)}m_{(1)}S^{-1}(h_{(2)})\ot h_{(1)})\\
&&\s(h_{(9)}S^{-1}(h_{(8)})\ot h_{(5)}m_{(2)})\s^{-1}(h_{(10)}\ot S^{-1}(h_{(7)})h_{(6)}m_{(3)})\\
&=&\S(h_{(3)}\cdot m_{(0)})\s(h_{(4)}m_{(1)}S^{-1}(h_{(2)})\ot h_{(1)})\s(h_{(5)}\ot m_{(2)})\\
&=&h\rhu m,
\end{array}$$
where we use (CQT$4'$) and (\ref{+-}) to obtain the fourth equality, use (\ref{lcocycle}) to obtain the sixth equality
and use (\ref{+-}) to get the seventh equality.
Again, from \leref{3.1}(b) one gets $m\in\hspace{0.05cm}_{\diamond}\u{\s}(M)$, and so $_{\diamond}M\subseteq\hspace{0.05cm}_{\diamond}\u{\s}(M)$.
Now replacing $H$, $\s$ and $M$ with $H^{\s}$, $\s^{-1}$ and $\u{\s}(M)$ respectively, one also gets $_{\diamond}\u{\s}(M)\subseteq\hspace{0.05cm}_{\diamond}\us^{-1}(\u{\s}(M))=\hspace{0.05cm}_{\diamond}M$. Hence we have $_{\diamond}\u{\s}(M)=\hspace{0.05cm}_{\diamond}M$.

If $A$ is a {\rm YD} $H$-module algebra, then $A_{\diamond}$ and $_{\diamond}A$ are subalgebras of $A$.
Thus the last statement follows immediately.
\end{proof}

Given two {\rm YD} $H$-modules $M$ and $N$, a generalized cotensor product  $M\wedge N$ was introduced in \cite{Z}. That is,
$$ M\wedge N=\left\{\S m_i\ot n_i\in M\ot N\left|\S(m_i\triangleleft\!\!-h)\ot n_i =\S m_i\ot(h-\!\!\triangleright n_i), \forall h\in{\mathcal H}_R\right.\right\}.$$
Observe that $M\wedge N$ is still an ${\mathcal H}_R$-bimodule with the left and right ${\mathcal H}_R$-module structures stemming from the left ${\mathcal H}_R$-module structure of $M$ and the right ${\mathcal H}_R$-module structure of $N$. In fact, this ${\mathcal H}_R$-bimodule structure of $M\wedge N$ comes from a {\rm YD} $H$-module structure on $M\wedge N$ by (\ref{LM})
and (\ref{RM}). The {\rm YD} $H$-module structure of $M\wedge N$ can be described as follows:

The left $H$-action on $M\wedge N$ is given by
\begin{equation}\label{wlm}
h\cdot\S(m_i\ot n_i)=\S h_{(1)}\cdot m_i\ot h_{(2)}\vartriangleright_1n_i
=\S h_{(1)}\vartriangleright_2m_i\ot h_{(2)}\cdot n_i,
\end{equation}
where $\S m_i\ot n_i\in M\wedge N$ and $h\in H$. The right $H$-comodule structure of $M\wedge N$
inherites from $M\ot N$. That is,
\begin{equation}\label{wrcom}
M\wedge N\ra(M\wedge N)\ot H,\hspace{0.1cm}\S m_i\ot n_i\mapsto
\S(m_{i(0)}\ot n_{i(0)})\ot n_{i(1)}m_{i(1)}.
\end{equation}

Let $\S m_i\ot n_i\in M\ot N$. Then by \cite[Lemma 2.9]{Z} we know that
$\S m_i\ot n_i\in M\wedge N$ if and only if
\begin{equation}\label{wedge}
\S h_{(1)}\cdot m_i\ot h_{(2)}\vartriangleright_1n_i
=\S h_{(1)}\vartriangleright_2m_i\ot h_{(2)}\cdot n_i,\hspace{0.1cm}\forall h\in H.
\end{equation}

Now we show that the equivalence functor $\us$ preserves the generalized cotensor product $\wedge$.

\begin{lemma}\lelabel{3.4}
Let $M$ and $N$ be {\rm YD} $H$-modules. Then $\u{\s}(M)\wedge\u{\s}(N)\cong\u{\s}(M\wedge N)$
as {\rm YD} $H^{\s}$-modules.
\end{lemma}

\begin{proof} By \thref{2.3} and its proof, we have a {\rm YD} $H^{\s}$-module isomorphism
$$\eta^{-1}_{M,N}: \u{\s}(M\ot N)\ra\u{\s}(M)\ot\u{\s}(N),\hspace{0.1cm}
x\ot y\mapsto\S x_{(0)}\ot y_{(0)}\s(y_{(1)}\ot x_{(1)}).$$
We show that $\eta^{-1}_{M,N}$ restricts to an isomorphism from $\u{\s}(M\wedge N)$ to $\u{\s}(M)\wedge\u{\s}(N)$. In order to simplify the computation,
we  write $x\ot y$ for an element $\S x_i\ot y_i$ in $M\wedge N$.
Let $x\ot y\in M\wedge N$. We show that the element $\eta^{-1}_{M,N}(x\ot y)=\S x_{(0)}\ot y_{(0)}\s(y_{(1)}\ot x_{(1)})$ sits in $\us(M)\wedge \us(N)$.
By the proof of \cite[Proposition 2.10]{Z},
$M\wedge N$ is an $H$-subcomdule of $M\ot N$. Hence
$\S(x_{(0)}\ot y_{(0)})\ot y_{(1)}x_{(1)}\in(M\wedge N)\ot H$. Following (\ref{wedge})
we have the identity for all $h\in H$:
\begin{equation}\label{w3.4}
\S h_{(1)}\cdot x_{(0)}\ot h_{(2)}\vartriangleright_1y_{(0)}\ot y_{(1)}x_{(1)}
=\S h_{(1)}\vartriangleright_2x_{(0)}\ot h_{(2)}\cdot y_{(0)}\ot y_{(1)}x_{(1)}.
\end{equation}
Now for any $h\in H^{\s}$, we verify the equation (\ref{wedge}) for the element $\eta^{-1}_{M,N}(x\ot y)$:
$$\begin{array}{cl}
&\S h_{(1)}\rhu x_{(0)}\ot h_{(2)}\vartriangleright_{1}^{\s}y_{(0)}\s(y_{(1)}\ot x_{(1)})\\
=&\S (h_{(2)}\cdot x_{(0)})_{(0)}\ot y_{(0)}\s((h_{(2)}\cdot x_{(0)})_{(1)}\ot h_{(1)})
\s^{-1}(h_{(3)}\ot x_{(1)})\\
&\s(y_{(1)}\ot h_{(4)})R(h_{(5)}\ot y_{(2)})\s^{-1}(h_{(6)}\ot y_{(3)})\s(y_{(4)}\ot x_{(2)})\\
\stackrel{(\ref{+-})}{=}&\S (h_{(2)}\cdot x_{(0)})_{(0)}\ot y_{(0)}
\s((h_{(2)}\cdot x_{(0)})_{(1)}\ot h_{(1)})\s^{-1}(h_{(3)}\ot x_{(1)})\\
&\s(y_{(1)}\ot h_{(4)})R(h_{(5)}\ot y_{(2)})\s(h_{(6)}y_{(3)}\ot x_{(2)})
\s^{-1}(h_{(7)}\ot y_{(4)}x_{(3)})\\
\stackrel{\rm (CQT4)}{=}&\S (h_{(2)}\cdot x_{(0)})_{(0)}\ot y_{(0)}
\s((h_{(2)}\cdot x_{(0)})_{(1)}\ot h_{(1)})\s^{-1}(h_{(3)}\ot x_{(1)})\\
&\s(y_{(1)}\ot h_{(4)})\s(y_{(2)}h_{(5)}\ot x_{(2)})R(h_{(6)}\ot y_{(3)})
\s^{-1}(h_{(7)}\ot y_{(4)}x_{(3)})\\
\stackrel{(\ref{lcocycle})}{=}&\S (h_{(2)}\cdot x_{(0)})_{(0)}\ot y_{(0)}
\s((h_{(2)}\cdot x_{(0)})_{(1)}\ot h_{(1)})\s^{-1}(h_{(3)}\ot x_{(1)})\\
&\s(h_{(4)}\ot x_{(2)})\s(y_{(1)}\ot h_{(5)}x_{(3)})R(h_{(6)}\ot y_{(2)})
\s^{-1}(h_{(7)}\ot y_{(3)}x_{(4)})\\
\stackrel{\rm(CQT4'')}{=}&\S (h_{(2)}\cdot x_{(0)})_{(0)}\ot y_{(0)}
\s((h_{(2)}\cdot x_{(0)})_{(1)}\ot h_{(1)})R(h_{(5)}\ot y_{(1)})\\
&\s(h_{(6)}y_{(2)}S^{-1}(h_{(4)})\ot h_{(3)}x_{(1)})
\s^{-1}(h_{(7)}\ot y_{(3)}x_{(2)})\\
=&\S (h_{(2)}\cdot x_{(0)})_{(0)}\ot(h_{(5)}\vartriangleright_1y_{(0)})
\s((h_{(2)}\cdot x_{(0)})_{(1)}\ot h_{(1)})\\
&\s(h_{(6)}y_{(1)}S^{-1}(h_{(4)})\ot h_{(3)}x_{(1)})
\s^{-1}(h_{(7)}\ot y_{(2)}x_{(2)})\\
\stackrel{(2)(3)}{=}&\S (h_{(3)}\cdot x_{(0)})_{(0)}\ot(h_{(4)}\vartriangleright_1y_{(0)})_{(0)}
\s((h_{(3)}\cdot x_{(0)})_{(1)}\ot h_{(1)})\\
&\s((h_{(4)}\vartriangleright_1y_{(0)})_{(1)}\ot(h_{(3)}\cdot x_{(0)})_{(2)}h_{(2)})
\s^{-1}(h_{(5)}\ot y_{(1)}x_{(1)})\\
\stackrel{(\ref{w3.4})}{=}&\S (h_{(3)}\vartriangleright_2x_{(0)})_{(0)}\ot(h_{(4)}\cdot y_{(0)})_{(0)}
\s((h_{(3)}\vartriangleright_2x_{(0)})_{(1)}\ot h_{(1)})\\
&\s((h_{(4)}\cdot y_{(0)})_{(1)}\ot(h_{(3)}\vartriangleright_2x_{(0)})_{(2)}h_{(2)})
\s^{-1}(h_{(5)}\ot y_{(1)}x_{(1)}).
\end{array}$$
On the other hand, we have
$$\begin{array}{cl}
&\S h_{(1)}\vartriangleright_{2}^{\s}x_{(0)}\ot h_{(2)}\rhu y_{(0)}\s(y_{(1)}\ot x_{(1)})\\
=&\S x_{(0)}\ot(h_{(3)}\cdot y_{(0)})_{(0)}R^{\s}(x_{(1)}\ot(S^{\s})^{-1}(h_{(1)}))\\
&\s((h_{(3)}\cdot y_{(0)})_{(1)}\ot h_{(2)})\s^{-1}(h_{(4)}\ot y_{(1)})\s(y_{(2)}\ot x_{(2)})\\

\stackrel{(\ref{+-})}{=}&\S x_{(0)}\ot(h_{(9)}\cdot y_{(0)})_{(0)}
\s(S^{-1}(h_{(2)})\ot h_{(1)})\s(S^{-1}(h_{(5)})\ot x_{(1)})\\
&R(x_{(2)}\ot S^{-1}(h_{(4)}))\s^{-1}(x_{(3)}\ot S^{-1}(h_{(3)}))\s^{-1}(h_{(7)}\ot S^{-1}(h_{(6)}))\\
&\s((h_{(9)}\cdot y_{(0)})_{(1)}\ot h_{(8)})
\s(h_{(10)}y_{(1)}\ot x_{(4)})\s^{-1}(h_{(11)}\ot y_{(2)}x_{(5)})\\
\stackrel{\rm(CQT4')(\ref{+-})}{=}&\S x_{(0)}\ot(h_{(12)}\cdot y_{(0)})_{(0)}
\s(x_{(3)}S^{-1}(h_{(4)})\ot h_{(1)})
\s^{-1}(x_{(4)}\ot S^{-1}(h_{(3)})h_{(2)})\\
&R(x_{(1)}\ot S^{-1}(h_{(6)}))\s(S^{-1}(h_{(8)})\ot h_{(7)}x_{(2)}S^{-1}(h_{(5)}))
\s^{-1}(h_{(10)}\ot S^{-1}(h_{(9)}))\\
&\s((h_{(12)}\cdot y_{(0)})_{(1)}\ot h_{(11)})
\s(h_{(13)}y_{(1)}\ot x_{(5)})\s^{-1}(h_{(14)}\ot y_{(2)}x_{(6)})\\
\stackrel{(2)}{=}&\S h_{(4)}\vartriangleright_2x_{(0)}\ot(h_{(11)}\cdot y_{(0)})_{(0)}
\s(x_{(2)}S^{-1}(h_{(2)})\ot h_{(1)})\\
&\s(S^{-1}(h_{(6)})\ot h_{(5)}x_{(1)}S^{-1}(h_{(3)}))
\s^{-1}(h_{(8)}\ot S^{-1}(h_{(7)}))\\
&\s((h_{(11)}\cdot y_{(0)})_{(1)}\ot h_{(9)})
\s((h_{(11)}\cdot y_{(0)})_{(2)}h_{(10)}\ot x_{(3)})\s^{-1}(h_{(12)}\ot y_{(1)}x_{(4)})\\
\stackrel{(3)(\ref{lcocycle})}{=}&\S (h_{(3)}\vartriangleright_2x_{(0)})_{(0)}\ot(h_{(9)}\cdot y_{(0)})_{(0)}
\s(x_{(1)}S^{-1}(h_{(2)})\ot h_{(1)})\\
&\s^{-1}(h_{(6)}\ot S^{-1}(h_{(5)}))
\s(S^{-1}(h_{(4)})\ot (h_{(3)}\vartriangleright_2x_{(0)})_{(1)})\\
&\s(h_{(7)}\ot x_{(2)})
\s((h_{(9)}\cdot y_{(0)})_{(1)}\ot h_{(8)}x_{(3)})
\s^{-1}(h_{(10)}\ot y_{(1)}x_{(4)})\\
\stackrel{(\ref{+-})}{=}&\S (h_{(3)}\vartriangleright_2x_{(0)})_{(0)}\ot(h_{(10)}\cdot y_{(0)})_{(0)}
\s(x_{(1)}S^{-1}(h_{(2)})\ot h_{(1)})\\
&\s(h_{(6)}S^{-1}(h_{(5)})\ot(h_{(3)}\vartriangleright_2x_{(0)})_{(1)})
\s^{-1}(h_{(7)}\ot S^{-1}(h_{(4)})(h_{(3)}\vartriangleright_2x_{(0)})_{(2)})\\
&\s(h_{(8)}\ot x_{(2)})
\s((h_{(10)}\cdot y_{(0)})_{(1)}\ot h_{(9)}x_{(3)})
\s^{-1}(h_{(11)}\ot y_{(1)}x_{(4)})\\
\end{array}$$
$$\begin{array}{cl}
=&\S (h_{(5)}\vartriangleright_2x_{(0)})_{(0)}\ot(h_{(10)}\cdot y_{(0)})_{(0)}
\s(x_{(1)}S^{-1}(h_{(2)})\ot h_{(1)})\\
&\s^{-1}(h_{(7)}\ot S^{-1}(h_{(6)})(h_{(5)}\vartriangleright_2x_{(0)})_{(1)}h_{(4)}S^{-1}(h_{(3)}))\\
&\s(h_{(8)}\ot x_{(2)})
\s((h_{(10)}\cdot y_{(0)})_{(1)}\ot h_{(9)}x_{(3)})
\s^{-1}(h_{(11)}\ot y_{(1)}x_{(4)})\\
\stackrel{(2)}{=}&\S (h_{(4)}\vartriangleright_2x_{(0)})\ot(h_{(8)}\cdot y_{(0)})_{(0)}\\
&\s^{-1}(h_{(5)}\ot x_{(1)}S^{-1}(h_{(3)}))\s(x_{(2)}S^{-1}(h_{(2)})\ot h_{(1)})\\
&\s(h_{(6)}\ot x_{(3)})
\s((h_{(8)}\cdot y_{(0)})_{(1)}\ot h_{(7)}x_{(4)})
\s^{-1}(h_{(9)}\ot y_{(1)}x_{(5)})\\
\stackrel{(\ref{+-})}{=}&\S (h_{(5)}\vartriangleright_2x_{(0)})\ot(h_{(10)}\cdot y_{(0)})_{(0)}\\
&\s(h_{(6)} x_{(1)}S^{-1}(h_{(4)})\ot h_{(1)})\s^{-1}(h_{(7)}\ot x_{(2)}S^{-1}(h_{(3)})h_{(2)})\\
&\s(h_{(8)}\ot x_{(3)})
\s((h_{(10)}\cdot y_{(0)})_{(1)}\ot h_{(9)}x_{(4)})
\s^{-1}(h_{(11)}\ot y_{(1)}x_{(5)})\\
\stackrel{(3)}{=}&\S (h_{(2)}\vartriangleright_2x_{(0)})_{(0)}\ot(h_{(4)}\cdot y_{(0)})_{(0)}
\s((h_{(2)}\vartriangleright_2x_{(0)})_{(1)}\ot h_{(1)})\\
&\s((h_{(4)}\cdot y_{(0)})_{(1)}\ot h_{(3)}x_{(1)})
\s^{-1}(h_{(5)}\ot y_{(1)}x_{(2)})\\
\stackrel{(2)}{=}&\S (h_{(3)}\vartriangleright_2x_{(0)})_{(0)}\ot(h_{(4)}\cdot y_{(0)})_{(0)}
\s((h_{(3)}\vartriangleright_2x_{(0)})_{(1)}\ot h_{(1)})\\
&\s((h_{(4)}\cdot y_{(0)})_{(1)}\ot (h_{(3)}\vartriangleright_2x_{(0)})_{(2)}h_{(2)})
\s^{-1}(h_{(5)}\ot y_{(1)}x_{(1)}).\\
\end{array}$$
Hence $\S h_{(1)}\rhu x_{(0)}\ot h_{(2)}\vartriangleright_{1}^{\s}y_{(0)}\s(y_{(1)}\ot x_{(1)}) =\S h_{(1)}\vartriangleright_{2}^{\s}x_{(0)}\ot h_{(2)}\rhu y_{(0)}\s(y_{(1)}\ot x_{(1)})$.
It follows from (\ref{wedge}) that $\eta^{-1}_{M,N}(x\ot y)\in\u{\s}(M)\wedge\u{\s}(N)$, and so $\eta^{-1}_{M,N}(M\wedge N)\subseteq\u{\s}(M)\wedge\u{\s}(N)$.

Now replacing $H$, $\s$, $M$ and $N$ with $H^{\s}$, $\s^{-1}$, $\u{\s}(M)$ and $\u{\s}(N)$ in the above,
respectively, one obtains $\eta_{M,N}(\u{\s}(M)\wedge\u{\s}(N))
\subseteq\us^{-1}(\u{\s}(M))\wedge\us^{-1}(\u{\s}(N))=M\wedge N$.
This implies $\u{\s}(M)\wedge\u{\s}(N)\subseteq \eta^{-1}_{M,N}(M\wedge N)$, and hence
$\eta^{-1}_{M,N}(\u{\s}(M\wedge N))=\eta^{-1}_{M,N}(M\wedge N)=\u{\s}(M)\wedge\u{\s}(N)$.
Thus the isomorphism $\eta^{-1}_{M,N}$ restricts to an isomorphism of $H^{\s}$-comodules
from $\u{\s}(M\wedge N)$ to $\u{\s}(M)\wedge\u{\s}(N)$
as $\eta^{-1}_{M,N}$ is obviously an $H^{\s}$-comodule homomorphism.

Finally, we verify that the restriction of $\eta^{-1}_{M,N}$ on $\us(M\wedge N)$ is a left $H^{\s}$-module map and consequently  a {\rm YD} $H^{\s}$-module isomorphism. Let $x\ot y\in\u{\s}(M\wedge N)$
and $h\in H^{\s}$. Then by \leref{2.1} and (\ref{wlm}) we have
$$\begin{array}{rcl}
h\cdot \eta^{-1}_{M,N}(x\ot y)&=&\S h_{(1)}\rhu x_{(0)}\ot h_{(2)}\vartriangleright_1^{\s}y_{(0)}
\s(y_{(1)}\ot x_{(1)})\\
&=&\S (h_{(2)}\cdot x_{(0)})_{(0)}\ot y_{(0)}\s((h_{(2)}\cdot x_{(0)})_{(1)}\ot h_{(1)})\\
&&\s^{-1}(h_{(3)}\ot x_{(1)})R^{\s}(h_{(4)}\ot y_{(1)})\s(y_{(2)}\ot x_{(2)})\\
&=&\S (h_{(2)}\cdot x_{(0)})_{(0)}\ot y_{(0)}\s((h_{(2)}\cdot x_{(0)})_{(1)}\ot h_{(1)})\\
&&\s^{-1}(h_{(3)}\ot x_{(1)})\s(y_{(1)}\ot h_{(4)})R(h_{(5)}\ot y_{(2)})\\
&&\s^{-1}(h_{(6)}\ot y_{(3)})\s(y_{(4)}\ot x_{(2)})\\

&\stackrel{(\ref{+-})}{=}&\S (h_{(2)}\cdot x_{(0)})_{(0)}\ot y_{(0)}
\s((h_{(2)}\cdot x_{(0)})_{(1)}\ot h_{(1)})\\
&&\s^{-1}(h_{(3)}\ot x_{(1)})\s(y_{(1)}\ot h_{(4)})R(h_{(5)}\ot y_{(2)})\\
&&\s(h_{(6)}y_{(3)}\ot x_{(2)})\s^{-1}(h_{(7)}\ot y_{(4)}x_{(3)})\\

&\stackrel{\rm(CQT4)}{=}&\S (h_{(2)}\cdot x_{(0)})_{(0)}\ot y_{(0)}
\s((h_{(2)}\cdot x_{(0)})_{(1)}\ot h_{(1)})\\
&&\s^{-1}(h_{(3)}\ot x_{(1)})\s(y_{(1)}\ot h_{(4)})
\s(y_{(2)}h_{(5)}\ot x_{(2)})\\
&&R(h_{(6)}\ot y_{(3)})\s^{-1}(h_{(7)}\ot y_{(4)}x_{(3)})\\
\end{array}$$
$$\begin{array}{rcl}
&\stackrel{(\ref{lcocycle})}{=}&\S (h_{(2)}\cdot x_{(0)})_{(0)}\ot y_{(0)}
\s((h_{(2)}\cdot x_{(0)})_{(1)}\ot h_{(1)})\\
&&\s^{-1}(h_{(3)}\ot x_{(1)})
\s(h_{(4)}\ot x_{(2)})\s(y_{(1)}\ot h_{(5)}x_{(3)})\\
&&R(h_{(6)}\ot y_{(2)})\s^{-1}(h_{(7)}\ot y_{(3)}x_{(4)})\\
&=&\S (h_{(2)}\cdot x_{(0)})_{(0)}\ot y_{(0)}
\s((h_{(2)}\cdot x_{(0)})_{(1)}\ot h_{(1)})\\
&&\s(y_{(1)}\ot h_{(3)}x_{(1)})
R(h_{(4)}\ot y_{(2)})\s^{-1}(h_{(5)}\ot y_{(3)}x_{(2)})\\

&\stackrel{(2)\rm(CQT4'')}{=}&\S (h_{(3)}\cdot x_{(0)})_{(0)}\ot y_{(0)}
\s((h_{(3)}\cdot x_{(0)})_{(1)}\ot h_{(1)})\\
&&R(h_{(5)}\ot y_{(1)})\s(h_{(6)}y_{(2)}S^{-1}(h_{(4)})\ot (h_{(3)}\cdot x_{(0)})_{(2)}h_{(2)})\\
&&\s^{-1}(h_{(7)}\ot y_{(3)}x_{(1)})\\

&=&\S (h_{(3)}\cdot x_{(0)})_{(0)}\ot h_{(5)}\vartriangleright_1y_{(0)}
\s((h_{(3)}\cdot x_{(0)})_{(1)}\ot h_{(1)})\\
&&\s(h_{(6)}y_{(1)}S^{-1}(h_{(4)})\ot (h_{(3)}\cdot x_{(0)})_{(2)}h_{(2)})
\s^{-1}(h_{(7)}\ot y_{(2)}x_{(1)})\\

&\stackrel{(3)}{=}&\S (h_{(3)}\cdot x_{(0)})_{(0)}\ot(h_{(4)}\vartriangleright_1y_{(0)})_{(0)}
\s((h_{(3)}\cdot x_{(0)})_{(1)}\ot h_{(1)})\\
&&\s((h_{(4)}\vartriangleright_1y_{(0)})_{(1)}\ot (h_{(3)}\cdot x_{(0)})_{(2)}h_{(2)})
\s^{-1}(h_{(5)}\ot y_{(1)}x_{(1)}).\\
\end{array}$$
On the other hand, we have
$$\begin{array}{rcl}
\eta^{-1}_{M,N}(h\rhu(x\ot y))&=&\S \eta^{-1}_{M,N}(h_{(3)}\cdot(x_{(0)}\ot y_{(0)}))
\s(h_{(4)}y_{(1)}x_{(1)}S^{-1}(h_{(2)})\ot h_{(1)})\\
&&\s^{-1}(h_{(5)}\ot y_{(2)}x_{(2)})\\
&=&\S (h_{(3)}\cdot x_{(0)})_{(0)}\ot(h_{(4)}\vartriangleright_1y_{(0)})_{(0)}\\
&&\s((h_{(4)}\vartriangleright_1y_{(0)})_{(1)}\ot (h_{(3)}\cdot x_{(0)})_{(1)})\\
&&\s(h_{(5)}y_{(1)}x_{(1)}S^{-1}(h_{(2)})\ot h_{(1)})\s^{-1}(h_{(6)}\ot y_{(2)}x_{(2)})\\
&=&\S (h_{(3)}\cdot x_{(0)})_{(0)}\ot(h_{(6)}\vartriangleright_1y_{(0)})_{(0)}\\
&&\s((h_{(6)}\vartriangleright_1y_{(0)})_{(1)}\ot (h_{(3)}\cdot x_{(0)})_{(1)})\\
&&\s(h_{(7)}y_{(1)}S^{-1}(h_{(5)})h_{(4)}x_{(1)}S^{-1}(h_{(2)})\ot h_{(1)})
\s^{-1}(h_{(8)}\ot y_{(2)}x_{(2)})\\
&\stackrel{(3)}{=}&\S (h_{(2)}\cdot x_{(0)})_{(0)}\ot(h_{(3)}\vartriangleright_1y_{(0)})_{(0)}\\
&&\s((h_{(3)}\vartriangleright_1y_{(0)})_{(1)}\ot (h_{(2)}\cdot x_{(0)})_{(1)})\\
&&\s((h_{(3)}\vartriangleright_1y_{(0)})_{(2)}(h_{(2)}\cdot x_{(0)})_{(2)}\ot h_{(1)})
\s^{-1}(h_{(4)}\ot y_{(1)}x_{(1)})\\
&\stackrel{(\ref{lcocycle})}{=}&\S (h_{(3)}\cdot x_{(0)})_{(0)}\ot(h_{(4)}\vartriangleright_1y_{(0)})_{(0)}
\s((h_{(3)}\cdot x_{(0)})_{(1)}\ot h_{(1)})\\
&&\s((h_{(4)}\vartriangleright_1y_{(0)})_{(1)}\ot (h_{(3)}\cdot x_{(0)})_{(2)}h_{(2)})
\s^{-1}(h_{(5)}\ot y_{(1)}x_{(1)}).\\
\end{array}$$
It follows that $\eta^{-1}_{M,N}(h\rhu(x\ot y))=h\cdot \eta^{-1}_{M,N}(x\ot y)$.   That is,  $\eta^{-1}_{M,N}$ is an $H^{\s}$-module homomorphism, and whence $\u{\s}(M)\wedge\u{\s}(N)\cong\u{\s}(M\wedge N)$ as {\rm YD} $H^{\s}$-modules.
\end{proof}

Now let $A$ and $B$ be {\rm YD} $H$-module algebras. Then $A$ and $B$ can be regarded as algebras
in ${\mathcal M}^H_R$ by forgetting the $H$-module structures of $A$ and $B$.
Then we can consider the braided product of $A$ and $B$ in ${\mathcal M}^H_R$.
Denote by $A\#_RB$ the braided product to differ from the braided product $A\#B$ in $\YDH$.
It follows from \cite[Proposition 2.11]{Z} that $A\wedge B$ is a subalgebra of $A\#_RB$
and the algebra $A\wedge B$ is a {\rm YD} $H$-module algebra with the $H$-structures given by (\ref{wlm}) and (\ref{wrcom}). We show that the equivalence functor $\us$ preserves the generalized cotensor product of algebras.

\begin{proposition}\prlabel{3.5}
Let $A$ and $B$ be {\rm YD} $H$-module algebras. Then $\u{\s}(A)\wedge\u{\s}(B)$ and $\u{\s}(A\wedge B)$ are isomorphic as {\rm YD} $H^{\s}$-module algebras.
\end{proposition}

\begin{proof} By \leref{3.4} and its proof, it is enough to show that the map
$$\eta^{-1}_{M,N}: \u{\s}(A\ot B)\ra\u{\s}(A)\ot\u{\s}(B), a\ot b\mapsto\S a_{(0)}\ot b_{(0)}\s(b_{(1)}\ot a_{(1)})$$
is an algebra homomorphism from $\u{\s}(A\#_RB)$ to $\u{\s}(A)\#_{R^{\s}}\u{\s}(B)$.
Obviously, $\eta^{-1}_{M,N}(1\#1)=1\#1$.
Now let $a$, $a'\in A$ and $b$, $b'\in B$.
Denote by $\cdot$ the product in $\u{\s}(A)\#_{R^{\s}}\u{\s}(B)$. Then we have
$$\begin{array}{ll}
&\eta^{-1}_{M,N}(a\#b)\cdot \eta^{-1}_{M,N}(a'\#b')\\
=&\S(a_{(0)}\#b_{(0)})\cdot(a'_{(0)}\#b'_{(0)})
\s(b_{(1)}\ot a_{(1)})\s(b'_{(1)}\ot a'_{(1)})\\
=&\S (a_{(0)}\bullet a'_{(0)}\#b_{(0)}\bullet b'_{(0)})R^{\s}(a'_{(1)}\ot b_{(1)})
\s(b_{(2)}\ot a_{(1)})\s(b'_{(1)}\ot a'_{(2)})\\
=&\S (a_{(0)}a'_{(0)}\#b_{(0)}b'_{(0)})\s^{-1}(a'_{(1)}\ot a_{(1)})\s^{-1}(b'_{(1)}\ot b_{(1)})\\
&R^{\s}(a'_{(2)}\ot b_{(2)})\s(b_{(3)}\ot a_{(2)})\s(b'_{(2)}\ot a'_{(3)})\\
=&\S (a_{(0)}a'_{(0)}\#b_{(0)}b'_{(0)})\s^{-1}(a'_{(1)}\ot a_{(1)})\s^{-1}(b'_{(1)}\ot b_{(1)})
\s(b_{(2)}\ot a'_{(2)})\\
&R(a'_{(3)}\ot b_{(3)})\s^{-1}(a'_{(4)}\ot b_{(4)})\s(b_{(5)}\ot a_{(2)})\s(b'_{(2)}\ot a'_{(5)})\\
\stackrel{(\ref{+-})}{=}&
\S (a_{(0)}a'_{(0)}\#b_{(0)}b'_{(0)})\s^{-1}(a'_{(1)}\ot a_{(1)})\s(b'_{(1)}b_{(1)}\ot a'_{(2)})
\s^{-1}(b'_{(2)}\ot b_{(2)}a'_{(3)})\\
&R(a'_{(4)}\ot b_{(3)})\s(a'_{(5)}b_{(4)}\ot a_{(2)})\s^{-1}(a'_{(6)}\ot b_{(5)}a_{(3)})
\s(b'_{(3)}\ot a'_{(7)})\\
\stackrel{\rm(CQT4)}{=}&
\S (a_{(0)}a'_{(0)}\#b_{(0)}b'_{(0)})\s^{-1}(a'_{(1)}\ot a_{(1)})\s(b'_{(1)}b_{(1)}\ot a'_{(2)})
R(a'_{(3)}\ot b_{(2)})\\
&\s^{-1}(b'_{(2)}\ot a'_{(4)}b_{(3)})\s(a'_{(5)}b_{(4)}\ot a_{(2)})
\s^{-1}(a'_{(6)}\ot b_{(5)}a_{(3)})\s(b'_{(3)}\ot a'_{(7)})\\
\stackrel{({\rm CQT4}')(\ref{+-})}{=}&
\S (a_{(0)}a'_{(0)}\#b_{(0)}b'_{(0)})\s^{-1}(a'_{(1)}\ot a_{(1)})R(a'_{(2)}\ot b_{(3)})\\
&\s(b'_{(1)}b_{(1)}\ot S(b_{(2)})a'_{(3)}b_{(4)})\s(b'_{(2)}a'_{(4)}b_{(5)}\ot a_{(2)})\\
&\s^{-1}(b'_{(3)}\ot a'_{(5)}b_{(6)}a_{(3)})
\s^{-1}(a'_{(6)}\ot b_{(7)}a_{(4)})\s(b'_{(4)}\ot a'_{(7)})\\
\stackrel{(\ref{rcocycle})}{=}&
\S (a_{(0)}a'_{(0)}\#b_{(0)}b'_{(0)})\s^{-1}(a'_{(1)}\ot a_{(1)})R(a'_{(2)}\ot b_{(5)})\\
&\s(b'_{(1)}b_{(1)}\ot S(b_{(4)})a'_{(3)}b_{(6)})\s(b'_{(2)}b_{(2)}S(b_{(3)})a'_{(4)}b_{(7)}\ot a_{(2)})\\
&\s^{-1}(b'_{(3)}a'_{(5)}\ot b_{(8)}a_{(3)})
\s^{-1}(b'_{(4)}\ot a'_{(6)})\s(b'_{(5)}\ot a'_{(7)})\\\stackrel{(\ref{lcocycle})}{=}&
\S (a_{(0)}a'_{(0)}\#b_{(0)}b'_{(0)})\s^{-1}(a'_{(1)}\ot a_{(1)})R(a'_{(2)}\ot b_{(4)})\\
&\s(S(b_{(3)})a'_{(3)}b_{(5)}\ot a_{(2)})\s(b'_{(1)}b_{(1)}\ot S(b_{(2)})a'_{(4)}b_{(6)}a_{(3)})\\
&\s^{-1}(b'_{(2)}a'_{(5)}\ot b_{(7)}a_{(4)})\\

\stackrel{\rm(CQT4)}{=}&\S (a_{(0)}a'_{(0)}\#b_{(0)}b'_{(0)})\s^{-1}(a'_{(1)}\ot a_{(1)})
\s(S(b_{(3)})b_{(4)}a'_{(2)}\ot a_{(2)})\\
&R(a'_{(3)}\ot b_{(5)})
\s(b'_{(1)}b_{(1)}\ot S(b_{(2)})a'_{(4)}b_{(6)}a_{(3)})\\
&\s^{-1}(b'_{(2)}a'_{(5)}\ot b_{(7)}a_{(4)})\\
\stackrel{\rm(CQT4)}{=}&
\S (a_{(0)}a'_{(0)}\#b_{(0)}b'_{(0)})
\s(b'_{(1)}b_{(1)}\ot S(b_{(2)})b_{(3)}a'_{(1)}a_{(1)})\\
&R(a'_{(2)}\ot b_{(4)})\s^{-1}(b'_{(2)}a'_{(3)}\ot b_{(5)}a_{(2)})\\

=&\S (a_{(0)}a'_{(0)}\#b_{(0)}b'_{(0)})
\s(b'_{(1)}b_{(1)}\ot a'_{(1)}a_{(1)})\\
&R(a'_{(2)}\ot b_{(2)})\s^{-1}(b'_{(2)}a'_{(3)}\ot b_{(3)}a_{(2)}).
\end{array}$$
On the other hand, we have
$$\begin{array}{ll}
&\eta^{-1}_{M,N}((a\#b)\bullet(a'\#b'))\\
=&\eta^{-1}_{M,N}(\S(a_{(0)}\#b_{(0)})(a'_{(0)}\#b'_{(0)})
\s^{-1}(b'_{(1)}a'_{(1)}\ot b_{(1)}a_{(1)}))\\
=&\S \eta^{-1}_{M,N}(a_{(0)}a'_{(0)}\#b_{(0)}b'_{(0)})R(a'_{(1)}\ot b_{(1)})
\s^{-1}(b'_{(1)}a'_{(2)}\ot b_{(2)}a_{(1)})\\
=&\S (a_{(0)}a'_{(0)}\#b_{(0)}b'_{(0)})\s(b'_{(1)}b_{(1)}\ot a'_{(1)}a_{(1)})
R(a'_{(2)}\ot b_{(2)})
\s^{-1}(b'_{(2)}a'_{(3)}\ot b_{(3)}a_{(2)}).\\
\end{array}$$
Hence $\eta^{-1}_{M,N}(a\#b)\cdot \eta^{-1}_{M,N}(a'\#b')=\eta^{-1}_{M,N}((a\#b)\bullet(a'\#b'))$. Thus, $\eta^{-1}_{M,N}$ is an algebra homomorphism from $\u{\s}(A\#_RB)$ to $\u{\s}(A)\#_{R^{\s}}\u{\s}(B)$.
\end{proof}

Let $A$ be a {\rm YD} $H$-module algebra. Recall from \cite[Definition 3.1]{Z} that the extension $A/A_{\diamond}$ is said to be a right ${\mathcal H}^*_R$-Galois extension if the  $k$-linear map
$$\b^r: A\ot_{A_{\diamond}}A\ra A\ot {\mathcal H}^*_R,\hspace{0.2cm}
\b^r(a\ot b)=\S a^{(0)}b\ot a^{(1)}$$
is an isomorphism. Similarly, the extension $A/_{\diamond}A$ is said to be left Galois
if the $k$-linear map
$$\b^l: A\ot_{_{\diamond}A}A\ra {\mathcal H}^*_R\ot A,\hspace{0.2cm}
\b^l(a\ot b)=\S b^{(-1)}\ot ab^{(0)}$$
is an isomorphism. If, in addition, the subalgebra $_{\diamond}A$ (or $A_{\diamond}$) is trivial and $A$ is faithfully flat over $k$, then $A$ is called a left (or right) ${\mathcal H}^*_R$-{\it Galois object}.
Denote by $\mathcal{E(H}_R)$ the category of {\rm YD} $H$-module algebras which are
${\mathcal H}^*_R$-bigalois objects. The morphisms in $\mathcal{E(H}_R)$ are
{\rm YD} $H$-module algebra homomorphisms (or equivalently, isomorphisms). If $A$ and $B$ are two objects of $\mathcal{E(H}_R)$, then $A\wedge B$ is an object of $\mathcal{E(H}_R)$ by \cite[Proposition 3.2]{Z}.

Let $H^*$ be the dual Hopf algebra of $H$. Then $H^*$ is a {\rm YD} $H$-module algebra with the
$H$-structures as follows: For $h^*, p\in H^*$ and $h\in H$, define
\begin{equation}\label{unit}
\begin{array}{rcl}
h\cdot p&=&\S p_{(1)}\langle p_{(2)}, h\rangle,\hspace{0.5cm}H\mbox{-action},\\
h^*\cdot p&=&\S h^*_{(2)}pS^{-1}(h^*_{(1)}),\hspace{0.5cm}H\mbox{-coaction}
\end{array}
\end{equation}
where we use $S$ for the antipodes of both $H$ and $H^*$ in order to simplify the notations
and we will do the same in the sequel. By \cite[Lemma 3.3]{Z}, $H^*$ is an object in $\mathcal{E(H}_R)$.
Denote by $I$ the object $H^*$ described above. It follows from \cite[Proposition 3.4]{Z} that $\mathcal{E(H}_R)$ is a monoidal category with the product $\wedge$ and the unit $I$.

Similarly, we have the monoidal category $\mathcal{E(H}^{\s}_{R^{\s}})$ for the CQT Hopf algebra $(H^{\s}, R^{\s})$. Our main task in the rest of this section is to show that the equivalence functor $\us$ restricts to an equivalence monoidal functor from the monoidal category $\mathcal{E(H)}$ to $\mathcal{E(H}^{\s}_{R^{\s}})$.
To this end, we first look at the unit $I^{\s}$ of $\mathcal{E(H}^{\s}_{R^{\s}})$ and show that $\us$ sends the unit  $I$ of $\mathcal{E(H)}$ to $I^{\s}$.

Let $(H^{\s})^*$ be the dual Hopf algebra of $H^{\s}$.
Then $I^{\s}=(H^{\s})^*$ as an algebra. The $H^{\s}$-structures of $I^{\s}$ is given by
\begin{equation}\label{units}
\begin{array}{rcl}
h\cdot_{\s}p&=&\S p_{<1>}\langle p_{<2>}, h\rangle,\hspace{0.5cm}H^{\s}\mbox{-action},\\
h^*\cdot_{\s}p&=&\S h^*_{<2>}p(S^{\s})^{-1}(h^*_{<1>}),\hspace{0.5cm}H^{\s}\mbox{-coaction}
\end{array}
\end{equation}
where $h^*, p\in(H^{\s})^*$ and $h\in H^{\s}$, and we use the sigama notation
$\S h^*_{<1>}\ot h^*_{<2>}$ for the comultiplication of an element $h^*\in (H^{\s})^*$ to differ from the comultiplication $\S h^*_{(1)}\ot h^*_{(2)}$ of $h^*\in H^*$.
Identify $(H\ot H)^*$ with $H^*\ot H^*$, we may assume that $\s$, $\s^{-1}\in H^*\ot H^*$ and may write $\s=\S \s^{(1)}\ot\s^{(2)}$ and $\s^{-1}=\S(\s^{-1})^{(1)}\ot(\s^{-1})^{(2)}$ in $H^*\ot H^*$.
Note that $(H^{\s})^*=H^*$ as algebras. Then for any $h^*\in H^*=(H^{\s})^*$ we have
$$\S h^*_{<1>}\ot h^*_{<2>}=\S\s(h^*_{(1)}\ot h^*_{(2)})\s^{-1}
=\S\s^{(1)}h^*_{(1)}(\s^{-1})^{(1)}\ot\s^{(2)}h^*_{(2)}(\s^{-1})^{(2)}.$$

Let $\chi$ be a $k$-linear map defined by
\begin{equation}\label{chi}
\chi: H\ra H,\hspace{0.1cm}\chi(h)=\S\s^{-1}(h_{(4)}\ot S^{-1}(h_{(3)})h_{(1)})h_{(2)},
\hspace{0.1cm}h\in H.
\end{equation}
Consequently, we have a $k$-linear map dual to $\chi$
\begin{equation}\label{chi*}
\chi^*: H^*\ra H^*,\hspace{0.1cm}\langle\chi^*(h^*), h\rangle=\langle h^*, \chi(h)\rangle,
\hspace{0.1cm} h^*\in H^*, h\in H.
\end{equation}

\begin{lemma}\lelabel{3.6}
$\chi$ is a $k$-linear isomorphism from $H$ to $H$ with the inverse given by
$$
\chi^{-1}(h)=\S\s^{-1}(S^{-1}(h_{(5)})\ot h_{(1)})\s(S^{-1}(h_{(4)})\ot h_{(3)})h_{(2)},
\hspace{0.1cm}h\in H.$$
Consequently, $\chi^*$ is a $k$-linear isomorphism as well.
\end{lemma}

\begin{proof}
Let $\lambda$ be the $k$-linear map given by $\lambda(h)=\S\s^{-1}(S^{-1}(h_{(5)})\ot h_{(1)})\s(S^{-1}(h_{(4)})\ot h_{(3)})h_{(2)}$,
$h\in H$. Then for any $h\in H$ we have
$$\begin{array}{rcl}
(\lambda\chi)(h)&=&\S \s^{-1}(h_{(4)}\ot S^{-1}(h_{(3)})h_{(1)})\lambda(h_{(2)})\\
&=&\S \s^{-1}(h_{(8)}\ot S^{-1}(h_{(7)})h_{(1)})\s^{-1}(S^{-1}(h_{(6)})\ot h_{(2)})\\
&&\s(S^{-1}(h_{(5)})\ot h_{(4)})h_{(3)}\\

&\stackrel{(\ref{rcocycle})}{=}&\S \s^{-1}(h_{(7)}S^{-1}(h_{(6)})\ot h_{(1)})\s^{-1}(h_{(8)}\ot S^{-1}(h_{(5)}))\\
&&\s(S^{-1}(h_{(4)})\ot h_{(3)})h_{(2)}\\

&\stackrel{(\ref{+-})}{=}&\S \s(h_{(6)}S^{-1}(h_{(5)})\ot h_{(2)})
\s^{-1}(h_{(7)}\ot S^{-1}(h_{(4)})h_{(3)})h_{(1)}\\
&=&h.
\end{array}$$
Similarly, one can check that $(\chi\lambda)(h)=h$ for any $h\in H$. It follows that
$\chi$ is an isomorphism with $\chi^{-1}=\lambda$, and so $\chi^*$ is also an isomorphism.
\end{proof}
We are ready to show that $\us(I)\cong I^{\s}$ in $\mathcal{E(H}^{\s}_{R^{\s})}$.

\begin{lemma}\lelabel{3.7}
$\u{\s}(I)$ and $I^{\s}$ are isomorphic {\rm YD} $H^{\s}$-module algebras.
\end{lemma}

\begin{proof}
By \leref{3.6}, $\chi^*: \u{\s}(H^*)\ra (H^{\s})^*$ is a $k$-linear isomorphism.
Hence it is enough to show that $\chi^*$ is a {\rm YD} $H^{\s}$-module algebra homomorphism
from $\u{\s}(I)$ to $I^{\s}$. Let $h^*\in(H^{\s})^*$, $p\in\u{\s}(H^*)$ and $h\in H^{\s}$. Then we have
$$\begin{array}{rcl}
\langle\chi^*(h^*\cdot p), h\rangle&=&\langle h^*\cdot p, \chi(h)\rangle\\
&=&\S\s^{-1}(h_{(4)}\ot S^{-1}(h_{(3)})h_{(1)})
\langle h^*_{(2)}pS^{-1}(h^*_{(1)}), h_{(2)}\rangle\\
&=&\S\s^{-1}(h_{(6)}\ot S^{-1}(h_{(5)})h_{(1)})\langle h^*_{(2)}, h_{(2)}\rangle
\langle p, h_{(3)}\rangle \langle S^{-1}(h^*_{(1)}), h_{(4)}\rangle\\
&=&\S\s^{-1}(h_{(6)}\ot S^{-1}(h_{(5)})h_{(1)})\langle h^*, S^{-1}(h_{(4)})h_{(2)}\rangle
\langle p, h_{(3)}\rangle\\
\end{array}$$
and
$$\begin{array}{rcl}
\langle h^*\cdot_{\s}\chi^*(p), h\rangle&=&\S\langle h^*_{<2>}\chi^*(p)(S^{\s})^{-1}(h^*_{<1>}), h\rangle\\

&=&\S\langle h^*_{<2>}, h_{(1)}\rangle
\langle\chi^*(p), h_{(2)}\rangle
\langle (S^{\s})^{-1}(h^*_{<1>}), h_{(3)}\rangle\\

&=&\S\langle p, \chi(h_{(2)})\rangle\langle h^*, (S^{\s})^{-1}(h_{(3)})\cdot_{\s}h_{(1)}\rangle\\
&=&\S\s^{-1}(h_{(5)}\ot S^{-1}(h_{(4)})h_{(2)})\langle p, h_{(3)}\rangle\s(S^{-1}(h_{(7)})\ot h_{(6)})\\
&&\s^{-1}(h_{(10)}\ot S^{-1}(h_{(9)})) \langle h^*, S^{-1}(h_{(8)})\cdot_{\s}h_{(1)}\rangle\\
&=&\S\s^{-1}(h_{(7)}\ot S^{-1}(h_{(6)})h_{(4)})\langle p, h_{(5)}\rangle\s(S^{-1}(h_{(9)})\ot h_{(8)})\\
&&\s^{-1}(h_{(14)}\ot S^{-1}(h_{(13)})) \s(S^{-1}(h_{(12)})\ot h_{(1)})\\
&&\langle h^*, S^{-1}(h_{(11)})h_{(2)}\rangle
\s^{-1}(S^{-1}(h_{(10)})\ot h_{(3)})\\

&\stackrel{(\ref{+-})}{=}&\S\s^{-1}(h_{(8)}\ot S^{-1}(h_{(7)})h_{(5)})\langle p, h_{(6)}\rangle\s(S^{-1}(h_{(10)})\ot h_{(9)})\\
&&\s(h_{(15)}S^{-1}(h_{(14)})\ot h_{(1)}) \s^{-1}(h_{(16)}\ot S^{-1}(h_{(13)})h_{(2)})\\
&&\langle h^*, S^{-1}(h_{(12)})h_{(3)}\rangle
\s^{-1}(S^{-1}(h_{(11)})\ot h_{(4)})\\

&=&\S\s^{-1}(S^{-1}(h_{(12)})\ot h_{(8)}S^{-1}(h_{(7)})h_{(3)})\s^{-1}(h_{(9)}\ot S^{-1}(h_{(6)})h_{(4)})\\
&&\s(S^{-1}(h_{(11)})\ot h_{(10)})\s^{-1}(h_{(15)}\ot S^{-1}(h_{(14)})h_{(1)})\\
&&\langle h^*, S^{-1}(h_{(13)})h_{(2)}\rangle \langle p, h_{(5)}\rangle\\

&\stackrel{(\ref{rcocycle})}{=}&
\S\s^{-1}(S^{-1}(h_{(11)})h_{(6)}\ot S^{-1}(h_{(5)})h_{(3)})\s^{-1}(S^{-1}(h_{(10)})\ot h_{(7)})\\
&&\s(S^{-1}(h_{(9)})\ot h_{(8)})\s^{-1}(h_{(14)}\ot S^{-1}(h_{(13)})h_{(1)})\\
&&\langle h^*, S^{-1}(h_{(12)})h_{(2)}\rangle \langle p, h_{(4)}\rangle\\

&=&\S\s^{-1}(h_{(6)}\ot S^{-1}(h_{(5)})h_{(1)})
\langle h^*, S^{-1}(h_{(4)})h_{(2)}\rangle \langle p, h_{(3)}\rangle.\\
\end{array}$$
So $\chi^*(h^*\cdot p)=h^*\cdot_{\s}\chi^*(p)$ and  $\chi^*$ is an $H^{\s}$-comodule
homomorphism from $\u{\s}(I)$ to $I^{\s}$.

Now we choose a pair of dual bases $\{h_1, h_2, \cdots, h_n\}$ in $H$ and
$\{h^*_1, h^*_2, \cdots, h^*_n\}$ in $H^*$ as $H$ is finitely generated projective. Then the comodule structures of $I$ and $I^{\s}$ are given by
$$I\ra I\ot H,\hspace{0.2cm}p\mapsto \S_{i=1}^n h^*_i\cdot p\ot h_i=\S p_{(0)}\ot p_{(1)}$$
and
$$I^{\s}\ra I^{\s}\ot H^{\s},\hspace{0.2cm}p\mapsto \S_{i=1}^n h^*_i\cdot_{\s} p\ot h_i
,$$
respectively. Thus for any $p\in\u{\s}(H^*)$ and $h$, $x\in H^{\s}$, by \leref{2.1} we have
$$\begin{array}{rl}
& \langle\chi^*(h\rhu p), x\rangle\\
=&\S \langle h_{(3)}\cdot p_{(0)}, \chi(x)\rangle
\s(h_{(4)}p_{(1)}S^{-1}(h_{(2)})\ot h_{(1)})\s^{-1}(h_{(5)}\ot p_{(2)})\\

=&\S \langle h^*_i\cdot p, \chi(x)h_{(3)}\rangle
\s(h_{(4)}h_{i(1)}S^{-1}(h_{(2)})\ot h_{(1)})\s^{-1}(h_{(5)}\ot h_{i(2)})\\

=&\S \langle h^*_{i(2)}pS^{-1}(h^*_{i(1)}), x_{(2)}h_{(3)}\rangle \s^{-1}(x_{(4)}\ot S^{-1}(x_{(3)})x_{(1)})\\
&\s(h_{(4)}h_{i(1)}S^{-1}(h_{(2)})\ot h_{(1)})\s^{-1}(h_{(5)}\ot h_{i(2)})\\

=&\S \langle h^*_i, S^{-1}(x_{(4)}h_{(5)})x_{(2)}h_{(3)}\rangle
\langle p, x_{(3)}h_{(4)}\rangle
\s^{-1}(x_{(6)}\ot S^{-1}(x_{(5)})x_{(1)})\\
&\s(h_{(6)}h_{i(1)}S^{-1}(h_{(2)})\ot h_{(1)})\s^{-1}(h_{(7)}\ot h_{i(2)})\\

=&\S\langle p, x_{(4)}h_{(5)}\rangle \s^{-1}(x_{(8)}\ot S^{-1}(x_{(7)})x_{(1)})\\
&\s(h_{(8)}S^{-1}(x_{(6)}h_{(7)})x_{(2)}h_{(3)}S^{-1}(h_{(2)})\ot h_{(1)})\\
&\s^{-1}(h_{(9)}\ot S^{-1}(x_{(5)}h_{(6)})x_{(3)}h_{(4)})\\
=&\S\langle p, x_{(4)}h_{(3)}\rangle \s^{-1}(x_{(8)}\ot S^{-1}(x_{(7)})x_{(1)})
\s(S^{-1}(x_{(6)})x_{(2)}\ot h_{(1)})\\
&\s^{-1}(h_{(5)}\ot S^{-1}(x_{(5)}h_{(4)})x_{(3)}h_{(2)})\\
\stackrel{(\ref{+-})}{=}&
\S\langle p, x_{(4)}h_{(4)}\rangle\s(x_{(8)}S^{-1}(x_{(7)})x_{(1)}\ot h_{(1)})\\
&\s^{-1}(x_{(9)}\ot S^{-1}(x_{(6)})x_{(2)}h_{(2)})
\s^{-1}(h_{(6)}\ot S^{-1}(x_{(5)}h_{(5)})x_{(3)}h_{(3)})\\
=&\S\s(x_{(1)}\ot h_{(1)})
\s^{-1}(x_{(7)}\ot h_{(7)}S^{-1}(h_{(6)})S^{-1}(x_{(6)})x_{(2)}h_{(2)})\\
&\s^{-1}(h_{(8)}\ot S^{-1}(h_{(5)})S^{-1}(x_{(5)})x_{(3)}h_{(3)})
\langle p, x_{(4)}h_{(4)}\rangle\\

\stackrel{(\ref{rcocycle})}{=}&
\S\s(x_{(1)}\ot h_{(1)})\s^{-1}(x_{(5)}h_{(5)}\ot S^{-1}(x_{(4)}h_{(4)})x_{(2)}h_{(2)})\\
&\s^{-1}(x_{(6)}\ot h_{(6)}) \langle p, x_{(3)}h_{(3)}\rangle\\
=&\S\s(x_{(1)}\ot h_{(1)})\langle p, \chi(x_{(2)}h_{(2)})\rangle \s^{-1}(x_{(3)}\ot h_{(3)})\\

=&\langle\chi^*(p), x\cdot_{\s}h\rangle\\
=&\langle h\cdot_{\s}\chi^*(p), x\rangle.
\end{array}$$
It follows that $\chi^*$ is an $H^{\s}$-module homomorphism from $\u{\s}(I)$ to $I^{\s}$.

Finally, we verify that $\chi^*$ is an algebra map. Let $p$, $q\in\u{\s}(I)$ and
$h\in H^{\s}$. Then
$$\begin{array}{rcl}
\langle \chi^*(p\bullet q), h\rangle
&=&\S\langle p_{(0)}q_{(0)}, \chi(h)\rangle\s^{-1}(q_{(1)}\ot p_{(1)})\\
&=&\S\langle (h^*_i\cdot p)(h^*_j\cdot q), h_{(2)}\rangle
\s^{-1}(h_{(4)}\ot S^{-1}(h_{(3)})h_{(1)})\s^{-1}(h_j\ot h_i)\\

&=&\S\langle h^*_i, S^{-1}(h_{(4)})h_{(2)}\rangle
\langle p, h_{(3)}\rangle
\langle h^*_j, S^{-1}(h_{(7)})h_{(5)}\rangle
\langle q, h_{(6)}\rangle\\
&&\s^{-1}(h_{(9)}\ot S^{-1}(h_{(8)})h_{(1)})\s^{-1}(h_j\ot h_i)\\

&=&\S\langle p, h_{(3)}\rangle\langle q, h_{(8)}\rangle
\s^{-1}(h_{(11)}\ot S^{-1}(h_{(10)})h_{(6)}S^{-1}(h_{(5)})h_{(1)})\\
&&\s^{-1}(S^{-1}(h_{(9)})h_{(7)}\ot S^{-1}(h_{(4)})h_{(2)})\\

&\stackrel{(\ref{rcocycle})}{=}&\S\langle p, h_{(2)}\rangle\langle q, h_{(6)}\rangle
\s^{-1}(h_{(9)}S^{-1}(h_{(8)})h_{(4)}\ot S^{-1}(h_{(3)})h_{(1)})\\
&&\s^{-1}(h_{(10)}\ot S^{-1}(h_{(7)})h_{(5)})\\

&=&\S\langle p, h_{(2)}\rangle\langle q, h_{(6)}\rangle
\s^{-1}(h_{(4)}\ot S^{-1}(h_{(3)})h_{(1)})\\
&&\s^{-1}(h_{(8)}\ot S^{-1}(h_{(7)})h_{(5)})\\

&=&\S\langle p, \chi(h_{(1)})\rangle\langle q, \chi(h_{(2)})\rangle\\
&=&\langle \chi^*(p)\chi^*(q), h\rangle.
\end{array}$$
Hence $\chi^*(p\bullet q)=\chi^*(p)\chi^*(q)$ for all $p,q\in \us(I)$. It is  clear that $\chi^*(\e)=\e$. Therefore, $\chi^*$ is an algebra map from $\u{\s}(I)$ to $I^{\s}$.
\end{proof}

Let $\omega: H\ot H\ra k$ be a $k$-linear map. Then $\omega$ induces two $k$-linear maps from $H$ to $H^*$. They are,
$$\omega_l: H\ra H^*,\hspace{0.2cm}\langle\omega_l(h), x\rangle=\omega(h\ot x),\hspace{0.2cm}h, x\in H$$
and
$$\omega_r: H\ra H^*,\hspace{0.2cm}\langle\omega_r(h), x\rangle=\omega(x\ot h),\hspace{0.2cm}h, x\in H.$$
For the 2-cocycle $\sigma$ and its inverse $\s^{-1}$, we have four
$k$-linear maps $\s_l$, $\s_r, (\s^{-1})_l$ and $(\s^{-1})_r$.
Note that $(\s^{-1})_l$ (or $(\s^{-1})_r$) coincides with the
convolution inverse $\s_l^{-1}$ (or $\s_r^{-1}$) of $\s_l$ (or
$\s_r$) in Hom$(H,H^*)$.

\begin{lemma}\lelabel{3.8}
Let $A$ be a {\rm YD} $H$-module algebra. Then the $k$-linear map
$$\phi: A\ot H^*\ra A\ot H^*,\hspace{0.1cm}a\ot h^*\mapsto
\S a_{(0)}\ot\s_l(a_{(1)})_{(2)}h^*S^{-1}(\s_l(a_{(1)})_{(1)})$$
is an isomorphism with the inverse given by
$$\phi^{-1}: A\ot H^*\ra A\ot H^*,\hspace{0.1cm}a\ot h^*\mapsto
\S a_{(0)}\ot\s^{-1}_l(a_{(1)})_{(2)}h^*S^{-1}(\s^{-1}_l(a_{(1)})_{(1)}),$$
where $H^*$ is the dual Hopf algebra of $H$.
\end{lemma}

\begin{proof} Identify $A\ot H^*$ with ${\rm Hom}(H, A)$. For $a\in A$, $h^*\in H^*$ and $h\in H$, we then have\\
$(\phi^{-1}\phi)(a\ot h^*)=\S a_{(0)}\ot\s^{-1}_l(a_{(1)})_{(2)}
\s_l(a_{(2)})_{(2)}h^*S^{-1}(\s_l(a_{(2)})_{(1)})S^{-1}(\s^{-1}_l(a_{(1)})_{(1)})$.
It leads to the following equations:
$$\begin{array}{rcl}
[(\phi^{-1}\phi)(a\ot h^*)](h)&=&\S a_{(0)}\langle\s^{-1}_l(a_{(1)})_{(2)}, h_{(1)}\rangle
\langle\s_l(a_{(2)})_{(2)}, h_{(2)}\rangle\langle h^*, h_{(3)}\rangle\\
&&\langle S^{-1}(\s_l(a_{(2)})_{(1)}), h_{(4)}\rangle
\langle S^{-1}(\s^{-1}_l(a_{(1)})_{(1)}), h_{(5)}\rangle\\
&=&\S a_{(0)}\langle\s^{-1}_l(a_{(1)}), S^{-1}(h_{(5)})h_{(1)}\rangle\\
&&\langle\s_l(a_{(2)}), S^{-1}(h_{(4)})h_{(2)}\rangle
\langle h^*, h_{(3)}\rangle\\
&=&\S a_{(0)}\s^{-1}(a_{(1)}\ot S^{-1}(h_{(5)})h_{(1)})\\
&&\s(a_{(2)}\ot S^{-1}(h_{(4)})h_{(2)})\langle h^*, h_{(3)}\rangle\\
&=&a\langle h^*, h\rangle=(a\ot h^*)(h).
\end{array}$$
It follows  that $\phi^{-1}\phi=id$. Similarly, one can show that $\phi\phi^{-1}=id$.
\end{proof}

\begin{lemma}\lelabel{3.9}
Let $A$ be a {\rm YD} $H$-module algebra. Then the $k$-linear map
$$\psi: H^*\ot A\ra H^*\ot A,\hspace{0.1cm}h^*\ot a\mapsto
\S \s_r(a_{(1)})_{(2)}h^*S^{-1}(\s_r(a_{(1)})_{(1)})\ot a_{(0)}$$
is an isomorphism with the inverse given by
$$\psi^{-1}: H^*\ot A\ra H^*\ot A,\hspace{0.1cm}h^*\ot a\mapsto
\S \s^{-1}_r(a_{(1)})_{(2)}h^*S^{-1}(\s^{-1}_r(a_{(1)})_{(1)})\ot a_{(0)}.$$
\end{lemma}

\begin{proof} Similar to the proof of \leref{3.8}.
\end{proof}
Now we are ready to show that $\us$ sends any object of $\mathcal{E(H}_R)$ into $\mathcal{E(H}^{\s}_{R^{\s}})$.
\begin{proposition}\prlabel{3.10}
Let $A$ be a {\rm YD} $H$-module algebra. Then\\
$\rm(a)$  $\u{\s}(A)/\u{\s}(A)_{\diamond}$ is a right $(\mathcal{H}^{\s}_{R^{\s}})^*$-Galois extension
if and only if $A/A_{\diamond}$ is a right\\
\mbox{\hspace{0.5cm}}$\mathcal{H}^*_R$-Galois extension.\\
$\rm(b)$  $\u{\s}(A)/_{\diamond}\u{\s}(A)$ is a left $(\mathcal{H}^{\s}_{R^{\s}})^*$-Galois extension
if and only if $A/_{\diamond}A$ is a left\\
\mbox{\hspace{0.5cm}}$\mathcal{H}^*_R$-Galois extension.
\end{proposition}

\begin{proof} By \cite[Eq.(12)]{Z}, the Galois map
$\b^r: A\ot_{A_{\diamond}}A\ra A\ot\mathcal{H}_R^*$ is given by
$$\b^r(a\ot b)=\S a_{[0](0)}b\ot a_{[1]}S^{-1}(R_r( a_{[0](1)})),\hspace{0.2cm}a, b\in A.$$
Similarly, the Galois map $\b_{\s}^r: \u{\s}(A)\ot_{\u{\s}(A)_{\diamond}}\u{\s}(A)
\ra\u{\s}(A)\ot(\mathcal{H}^{\s}_{R^{\s}})^*$ is given by
$$\begin{array}{rcl}
\b_{\s}^r(a\ot b)&=&\S a_{<0>(0)}\bullet b\ot a_{<1>}(S^{\s})^{-1}(R^{\s}_r( a_{<0>(1)}))\\
&=&\S a_{<0>(0)}b_{(0)}\s^{-1}(b_{(1)}\ot a_{<0>(1)})\ot a_{<1>}(S^{\s})^{-1}(R^{\s}_r( a_{<0>(2)})).\\
\end{array}$$
By \leref{3.2} and \leref{3.3}, it is straightforward to check that the $k$-linear isomorphism
$A\ot A\ra\u{\s}(A)\ot\u{\s}(A),\hspace{0.1cm}a\ot b\mapsto\S a_{(0)}\ot b_{(0)}\s(b_{(1)}\ot a_{(1)})$
induces an isomorphism
$$f: A\ot_{A_{\diamond}}A\ra\u{\s}(A)\ot_{\u{\s}(A)_{\diamond}}\u{\s}(A),\hspace{0.1cm}
a\ot b\mapsto\S a_{(0)}\ot b_{(0)}\s(b_{(1)}\ot a_{(1)}).$$
As $k$-modules, we may regard
$\u{\s}(A)\ot(\mathcal{H}^{\s}_{R^{\s}})^*=A\ot\mathcal{H}_R^*=A\ot H^*$.
We claim that the following diagram is commutative:
$$\begin{array}{ccc}
A\ot_{A_{\diamond}}A&\stackrel{\b^r}{\longrightarrow}&A\ot\mathcal{H}_R^*\\
f\downarrow&&\hspace{1.5cm}\downarrow({\rm id}\ot\chi^*)\phi\\
\u{\s}(A)\ot_{\u{\s}(A)_{\diamond}}\u{\s}(A)&\stackrel{\b_{\s}^r}{\longrightarrow}&
\u{\s}(A)\ot(\mathcal{H}^{\s}_{R^{\s}})^*,
\end{array}
$$
where $\chi^*$ and $\phi$ are isomorphisms given in (\ref{chi*}) and \leref{3.8}, respectively. In fact, let $a, b\in A$. Then we have
$$\begin{array}{rcl}
(\b_{\s}^rf)(a\ot b)&=&\S a_{(0)<0>(0)}b_{(0)}\ot a_{(0)<1>}(S^{\s})^{-1}(R^{\s}_r( a_{(0)<0>(2)}))\\
&&\s^{-1}(b_{(1)}\ot a_{(0)<0>(1)})\s(b_{(2)}\ot a_{(1)})\\
\end{array}$$
and
\begin{eqnarray*}
(({\rm id}\ot\chi^*)\phi\b^r)(a\ot b)
&=& \S a_{[0](0)}b_{(0)}\ot\chi^*(\s_l(b_{(1)}a_{[0](1)})_{(2)}
a_{[1]}S^{-1}(R_r( a_{[0](2)})) \\
&& S^{-1}(\s_l(b_{(1)}a_{[0](1)})_{(1)})).
\end{eqnarray*}
Identify $A\ot H^*$ with ${\rm Hom}(H, A)$,  then we have for any $h\in H$,
$$\begin{array}{cl}
&[(\b_{\s}^rf)(a\ot b)](h)\\
=&\S a_{(0)<0>(0)}b_{(0)}\langle a_{(0)<1>}, h_{(1)}\rangle
\langle R^{\s}_r( a_{(0)<0>(2)}),(S^{\s})^{-1}(h_{(2)})\rangle\\
&\s^{-1}(b_{(1)}\ot a_{(0)<0>(1)})\s(b_{(2)}\ot a_{(1)})\\

=&\S (h_{(1)}\rhu a_{(0)})_{(0)}b_{(0)}\s^{-1}(b_{(1)}\ot(h_{(1)}\rhu a_{(0)})_{(1)})\\
&R^{\s}((S^{\s})^{-1}(h_{(2)})\ot(h_{(1)}\rhu a_{(0)})_{(2)})\s(b_{(2)}\ot a_{(1)})\\

=&\S (h_{(2)}\cdot a_{(0)})_{(0)}b_{(0)}\s^{-1}(b_{(1)}\ot(h_{(2)}\cdot a_{(0)})_{(1)})\\
&R^{\s}(S^{-1}(h_{(6)})\ot(h_{(2)}\cdot a_{(0)})_{(2)})\s(b_{(2)}\ot a_{(2)})
\s((h_{(2)}\cdot a_{(0)})_{(3)}\ot h_{(1)})\\
&\s^{-1}(h_{(3)}\ot a_{(1)})\s(S^{-1}(h_{(5)})\ot h_{(4)})\s^{-1}(h_{(8)}\ot S^{-1}(h_{(7)}))\\

=&\S (h_{(2)}\cdot a_{(0)})_{(0)}b_{(0)}\s^{-1}(b_{(1)}\ot(h_{(2)}\cdot a_{(0)})_{(1)})
\s((h_{(2)}\cdot a_{(0)})_{(2)}\ot S^{-1}(h_{(8)}))\\
&R(S^{-1}(h_{(7)})\ot(h_{(2)}\cdot a_{(0)})_{(3)})
\s^{-1}(S^{-1}(h_{(6)})\ot(h_{(2)}\cdot a_{(0)})_{(4)})\\
&\s((h_{(2)}\cdot a_{(0)})_{(5)}\ot h_{(1)})\s^{-1}(h_{(3)}\ot a_{(1)})\\
&\s(S^{-1}(h_{(5)})\ot h_{(4)})\s^{-1}(h_{(10)}\ot S^{-1}(h_{(9)}))\s(b_{(2)}\ot a_{(2)})\\

\stackrel{(\ref{+-})}{=}
&\S (h_{(3)}\cdot a_{(0)})_{(0)}b_{(0)}\s(b_{(1)}(h_{(3)}\cdot a_{(0)})_{(1)}\ot S^{-1}(h_{(11)}))\\
&\s^{-1}(b_{(2)}\ot (h_{(3)}\cdot a_{(0)})_{(2)}S^{-1}(h_{(10)}))
R(S^{-1}(h_{(9)})\ot(h_{(3)}\cdot a_{(0)})_{(3)})\\
&\s(S^{-1}(h_{(8)})(h_{(3)}\cdot a_{(0)})_{(4)}\ot h_{(1)})
\s^{-1}(S^{-1}(h_{(7)})\ot (h_{(3)}\cdot a_{(0)})_{(5)}h_{(2)})\\
&\s^{-1}(h_{(4)}\ot a_{(1)})\s(S^{-1}(h_{(6)})\ot h_{(5)})
\s^{-1}(h_{(13)}\ot S^{-1}(h_{(12)}))\s(b_{(3)}\ot a_{(2)})\\

\stackrel{\rm(2)(CQT4)}{=}
&\S (h_{(2)}\cdot a_{(0)})_{(0)}b_{(0)}\s(b_{(1)}(h_{(2)}\cdot a_{(0)})_{(1)}\ot S^{-1}(h_{(11)}))\\
&\s^{-1}(b_{(2)}\ot (h_{(2)}\cdot a_{(0)})_{(2)}S^{-1}(h_{(10)}))
\s((h_{(2)}\cdot a_{(0)})_{(3)}S^{-1}(h_{(9)})\ot h_{(1)})\\
&R(S^{-1}(h_{(8)})\ot(h_{(2)}\cdot a_{(0)})_{(4)})
\s^{-1}(S^{-1}(h_{(7)})\ot h_{(3)}a_{(1)})\\
&\s^{-1}(h_{(4)}\ot a_{(2)})\s(S^{-1}(h_{(6)})\ot h_{(5)})
\s^{-1}(h_{(13)}\ot S^{-1}(h_{(12)}))\s(b_{(3)}\ot a_{(3)})\\
\stackrel{(\ref{+-})(\ref{rcocycle})}{=}
&\S (h_{(3)}\cdot a_{(0)})_{(0)}b_{(0)}\s(b_{(1)}(h_{(3)}\cdot a_{(0)})_{(1)}\ot S^{-1}(h_{(13)}))\\
&\s(b_{(2)}(h_{(3)}\cdot a_{(0)})_{(2)}S^{-1}(h_{(12)})\ot h_{(1)})
\s^{-1}(b_{(3)}\ot (h_{(3)}\cdot a_{(0)})_{(3)}S^{-1}(h_{(11)})h_{(2)})\\

&R(S^{-1}(h_{(10)})\ot(h_{(3)}\cdot a_{(0)})_{(4)})
\s^{-1}(S^{-1}(h_{(9)})h_{(4)}\ot a_{(1)})\\
&\s^{-1}(S^{-1}(h_{(8)})\ot h_{(5)})
\s(S^{-1}(h_{(7)})\ot h_{(6)})
\s^{-1}(h_{(15)}\ot S^{-1}(h_{(14)}))\s(b_{(4)}\ot a_{(2)})\\
\stackrel{(\ref{lcocycle})({\rm CQT4})}{=}
&\S (h_{(4)}\cdot a_{(0)})_{(0)}b_{(0)}
\s^{-1}(h_{(10)}\ot S^{-1}(h_{(9)}))\s(S^{-1}(h_{(8)})\ot h_{(1)})\\
&\s(b_{(1)}(h_{(4)}\cdot a_{(0)})_{(1)}\ot S^{-1}(h_{(7)})h_{(2)})
R(S^{-1}(h_{(6)})\ot(h_{(4)}\cdot a_{(0)})_{(2)})\\
&\s^{-1}(b_{(2)}\ot S^{-1}(h_{(5)})(h_{(4)}\cdot a_{(0)})_{(3)}h_{(3)})
\s(b_{(3)}\ot a_{(1)})\\

\stackrel{(2)(\ref{+-})}{=}
&\S (h_{(4)}\cdot a_{(0)})_{(0)}b_{(0)}
\s(h_{(11)}S^{-1}(h_{(10)})\ot h_{(1)})\s^{-1}(h_{(12)}\ot S^{-1}(h_{(9)})h_{(2)})\\
&\s(b_{(1)}(h_{(4)}\cdot a_{(0)})_{(1)}\ot S^{-1}(h_{(8)})h_{(3)})
R(S^{-1}(h_{(7)})\ot(h_{(4)}\cdot a_{(0)})_{(2)})\\
&\s^{-1}(b_{(2)}\ot S^{-1}(h_{(6)})h_{(5)}a_{(1)})
\s(b_{(3)}\ot a_{(2)})\\
=&\S (h_{(3)}\cdot a)_{(0)}b_{(0)}
\s^{-1}(h_{(7)}\ot S^{-1}(h_{(6)})h_{(1)})\\
&\s(b_{(1)}(h_{(3)}\cdot a)_{(1)}\ot S^{-1}(h_{(5)})h_{(2)})
R(S^{-1}(h_{(4)})\ot(h_{(3)}\cdot a)_{(2)}),\\
\end{array}$$
and
$$\begin{array}{cl}
&[(({\rm id}\ot\chi^*)\phi\b^r)(a\ot b)](h)\\
=&\S a_{[0](0)}b_{(0)}
\langle\s_l(b_{(1)}a_{[0](1)})_{(2)}a_{[1]}S^{-1}(R_r( a_{[0](2)}))
S^{-1}(\s_l(b_{(1)}a_{[0](1)})_{(1)}), \chi(h)\rangle\\

=&\S a_{[0](0)}b_{(0)}\s^{-1}(h_{(7)}\ot S^{-1}(h_{(6)})h_{(1)})
\langle\s_l(b_{(1)}a_{[0](1)})_{(2)}, h_{(2)}\rangle\\
&\langle a_{[1]}, h_{(3)}\rangle\langle R_r( a_{[0](2)}), S^{-1}(h_{(4)})\rangle
\langle \s_l(b_{(1)}a_{[0](1)})_{(1)}, S^{-1}(h_{(5)})\rangle\\

=&\S a_{[0](0)}b_{(0)}\s^{-1}(h_{(7)}\ot S^{-1}(h_{(6)})h_{(1)})
\s(b_{(1)}a_{[0](1)}\ot S^{-1}(h_{(5)})h_{(2)})\\
&\langle a_{[1]}, h_{(3)}\rangle R(S^{-1}(h_{(4)})\ot a_{[0](2)})\\

=&\S (h_{(3)}\cdot a)_{(0)}b_{(0)}\s^{-1}(h_{(7)}\ot S^{-1}(h_{(6)})h_{(1)})\\
&\s(b_{(1)}(h_{(3)}\cdot a)_{(1)}\ot S^{-1}(h_{(5)})h_{(2)})
R(S^{-1}(h_{(4)})\ot (h_{(3)}\cdot a)_{(2)}).\\
\end{array}$$
Since $\chi^*$ and $\phi$ are isomorphisms, it follows from the foregoing commutative diagram that $\b^r_{\s}$ is an isomorphism if and only if $\b^r$ is an isomorphism.
This completes the proof of  Part (a).

Part (b) follows from a similar argument by using the isomorphisms $\psi$ in \leref{3.9} and $\chi^*$.
\end{proof}

Combining \leref{3.3}, \prref{3.5}, \leref{3.7} and \prref{3.10},
we now arrive at the monoidal equivalence between $\mathcal{E(H)}_R$ and $\mathcal{E(H}^{\s}_{R^{\s}})$.

\begin{theorem}\thlabel{3.11}
$\mathcal{E(H}^{\s}_{R^{\s}})$ and $\mathcal{E(H}_R)$ are equivalent monoidal categories.
\end{theorem}

Let $E(\mathcal{H}_R)$ be the set of the isomorphism classes of objects in $\mathcal{E(H}_R)$.
Then $E(\mathcal{H}_R)$ is a semigroup with the product induced by the generalized cotensor product $\wedge$.
Denote by ${\rm Gal}(\mathcal{H}_R)$ the subset of $E(\mathcal{H}_R)$ consisting of
the isomorphism classes of objects in $\mathcal{E(H}_R)$ that are quantum commutative
in $_H\mathcal{YD}^H$, where a YD $H$-module algebra $A$ is quantum commutative in case it satisfies
\begin{equation}\label{qc}
ab=\sum b_{(0)}(b_{(1)}\cdot a),\ \ \mathrm{for\ all}\ a, b\in A.
\end{equation}
The subset ${\rm Gal}(\mathcal{H}_R)$ is a group (see \cite[Theorem 3.9]{Z}). Since the quantum commutativity (\ref{qc}) is defined by the braiding of the category $\YDH$ and $\us$ is a braided monoidal functor, $\us(A)$ is quantum commutative in $\YDHs$ if $A$ is quantum commutative in $\YDH$. Thus we obtain the following main result of this section.

\begin{theorem}\thlabel{3.12} The functor $\us$ induces a group isomorphism from
 ${\rm Gal}(\mathcal{H}_R)$ to ${\rm Gal}(\mathcal{H}^{\s}_{R^{\s}})$ sending an element $[A]$ to the element $[\us(A)]$.
\end{theorem}

If $\s$ is a lazy  $2$-cocycle on $H$, then $\s$ induces a group isomorphism
from ${\rm Gal}(\mathcal{H}_R)$ to ${\rm Gal}(\mathcal{H}_{R^{\s}})$. For example, every CQT structure $R_t$ of $H_4$ in \exref{2.12} is equal to $R_0^{\s}$ for some lazy cocycle $\s$. By \thref{3.12}, we have $\Gal(\mathcal{H}_{R_t})\cong \Gal(\mathcal{H}_{R_0})$. Similar to \coref{2.7} we have that $\mathrm{H}^2_L(H)$ acts on ${\rm Gal}(\mathcal{H}_R)$ by automorphisms.

To end this section, we show that the exact sequence (\ref{exactsqn}) is stable under the equivalence functor $\us$.  For simplifying notations we will write in the sequel $M_0$ for the coinvariant subset $M^{coH}$ of right $H$-comodule $M$:
$$\left\{ m\in M \left| \rho(m)=\S m_{(0)}\ot m_{(1)}=m\ot 1\right.\right\}.$$
Let $A$ be a right $H^{op}$-comodule algebra. Recall that $A/A_0$ is called
{\it an $H^{op}$-Galois extension} if the Galois map is bijective:
$$\b: A\ot_{A_0}A\ra A\ot H^{op},\hspace{0.2cm}a\ot b\mapsto\S ab_{(0)}\ot b_{(1)}.$$
If, in addition, $A$ is $R$-Azumaya, then we say that $A$ is {\it a Galois $R$-Azumaya
algebra}. It is known that any element of ${\rm BC}(k, H, R)$
can be represented by a Galois $R$-Azumaya algebra \cite[Corollary 4.2]{Z}.

Now let $A$ be a right $H^{op}$-comodule algebra such that $A/A_0$ is Galois.
Denote by $\pi(A)$ the centralizer subalgebra $C_A(A_0)$ of $A_0$ in $A$.
It is clear that $\pi(A)$ is an $H^{op}$-subcomodule subalgebra of $A$.
The Miyashita-Ulbrich-Van Oystaeyen action \cite{Miya, Ul, Vo} of $H$ on $\pi(A)$ is given by
\begin{equation}\label{MiU}
h\cdot a=\S X_i(h)aY_i(h), \hspace{0.2cm}a\in\pi(A), h\in H,
\end{equation}
where $\S X_i(h)\ot Y_i(h)=\b^{-1}(1\ot h)$ for $h\in H$. It is well known
that the right $H^{op}$-comodule algebra $\pi(A)$ together
with the action (\ref{MiU}) is a {\rm YD} $H$-module algebra (e.g., see \cite{CVZ1,Ul}). Moreover, $\pi(A)$ is quantum commutative in $_H\mathcal{YD}^H$. If, in addition, $A$ is $R$-Azumaya, $\pi(A)$ is an object in ${\rm Gal}(\mathcal{H}_R)$ (see \cite[Proposition 4.6]{Z}). Thus,  $\pi$ induces a group homomorphism
$$\stackrel{\sim}{\pi}: {\rm BC}(k, H, R)\ra{\rm Gal}(\mathcal{H}_R),
\hspace{0.2cm}[A]\mapsto [\pi(A)],$$
where $[A]$ is taken in ${\rm BC}(k, H, R)$ such that $A$ is a Galois $R$-Azumaya algebra. The group homomorphism $\stackrel{\sim}{\pi}$ fits in the exact sequence (\ref{exactsqn}) of groups:
$$1\ra{\rm Br}(k)\ra{\rm BC}(k, H, R)\stackrel{\stackrel{\sim}{\pi}}{\ra}{\rm Gal}(\mathcal{H}_R).$$
Similarly, we have a group homomorphism
$\stackrel{\sim}{\pi}: {\rm BC}(k,H^{\s},R^{\s})\ra{\rm Gal}(\mathcal{H}^{\s}_{R^{\s}})$
and such an exact group sequence (\ref{exactsqn}) for ${\rm BC}(k,H^{\s},R^{\s})$ and
${\rm Gal}(\mathcal{H}^{\s}_{R^{\s}})$.
We will prove that the exact sequence (\ref{exactsqn}) is stable under a cocycle deformation. We will need the following isomorphism $\xi$ later on.

\begin{lemma}\lelabel{3.13}
Let $A$ be a right $H^{op}$-comodule algebra. Then the $k$-linear map
$$\xi: A\ot H\ra A\ot H,\hspace{0.2cm}a\ot h\mapsto
\S a_{(0)}\ot h_{(4)}\s(S^{-1}(h_{(2)})\ot h_{(1)})\s^{-1}(S^{-1}(h_{(3)})\ot a_{(1)})$$
is an isomorphism with the inverse given by
$$\xi^{-1}(a\ot h)=\S a_{(0)}\ot h_{(3)}\s^{-1}(h_{(2)}\ot S^{-1}(h_{(1)})a_{(1)}),
\hspace{0.2cm}a\in A, h\in H.$$
\end{lemma}

\begin{proof} Let $\varphi: A\ot H\ra A\ot H$ be a $k$-linear map given by
$$\varphi(a\ot h)=\S a_{(0)}\ot h_{(3)}\s^{-1}(h_{(2)}\ot S^{-1}(h_{(1)})a_{(1)}).$$
Then for any $a\ot h\in A\ot H$, we have
$$\begin{array}{rcl}
(\varphi\xi)(a\ot h)&=&\S\varphi(a_{(0)}\ot h_{(4)})
\s(S^{-1}(h_{(2)})\ot h_{(1)})\s^{-1}(S^{-1}(h_{(3)})\ot a_{(1)})\\
&=&\S (a_{(0)}\ot h_{(6)})\s(S^{-1}(h_{(2)})\ot h_{(1)})\\
&&\s^{-1}(h_{(5)}\ot S^{-1}(h_{(4)})a_{(1)})\s^{-1}(S^{-1}(h_{(3)})\ot a_{(2)})\\
&\stackrel{(\ref{rcocycle})}{=}&\S (a_{(0)}\ot h_{(7)})\s(S^{-1}(h_{(2)})\ot h_{(1)})\\
&&\s^{-1}(h_{(5)}S^{-1}(h_{(4)})\ot a_{(1)})\s^{-1}(h_{(6)}\ot S^{-1}(h_{(3)}))\\
&=&\S (a\ot h_{(5)})\s^{-1}(h_{(4)}\ot S^{-1}(h_{(3)}))\s(S^{-1}(h_{(2)})\ot h_{(1)})\\
&\stackrel{(\ref{+-})}{=}&\S (a\ot h_{(7)})
\s(h_{(5)}S^{-1}(h_{(4)})\ot h_{(1)})\s^{-1}(h_{(6)}\ot S^{-1}(h_{(3)})h_{(2)})\\
&=&a\ot h
\end{array}$$
and
$$\begin{array}{rcl}
(\xi\varphi)(a\ot h)&=&\S\xi(a_{(0)}\ot h_{(3)})\s^{-1}(h_{(2)}\ot S^{-1}(h_{(1)})a_{(1)})\\
&=&\S(a_{(0)}\ot h_{(6)})\s(S^{-1}(h_{(4)})\ot h_{(3)})\s^{-1}(S^{-1}(h_{(5)})\ot a_{(1)})\\
&&\s^{-1}(h_{(2)}\ot S^{-1}(h_{(1)})a_{(2)})\\
&=&\S(a_{(0)}\ot h_{(8)})\s(S^{-1}(h_{(6)})\ot h_{(5)})\s^{-1}(S^{-1}(h_{(7)})\ot a_{(1)})\\
&&\s(h_{(3)}S^{-1}(h_{(2)})\ot a_{(2)})\s^{-1}(h_{(4)}\ot S^{-1}(h_{(1)})a_{(3)})\\
&\stackrel{(\ref{+-})}{=}&\S(a_{(0)}\ot h_{(7)})\s(S^{-1}(h_{(5)})\ot h_{(4)})
\s^{-1}(S^{-1}(h_{(6)})\ot a_{(1)})\\
&&\s^{-1}(h_{(3)}\ot S^{-1}(h_{(2)}))\s(S^{-1}(h_{(1)})\ot a_{(2)})\\
&\stackrel{(\ref{+-e})}{=}&\S(a_{(0)}\ot h_{(3)})
\s^{-1}(S^{-1}(h_{(2)})\ot a_{(1)})\s(S^{-1}(h_{(1)})\ot a_{(2)})\\
&&=a\ot h.\\
\end{array}$$
This shows that $\xi$ is an isomorphism and $\xi^{-1}=\varphi$.
\end{proof}

Let $A$ be a right $H^{op}$-comodule algebra. Recall that the
product of $\us(A)$ is given by $a \bullet b=\sum a_{(0)}b_{(0)}
\s^{-1}(a_{(1)}\ot b_{(1)})$, for all $a, b\in A$. Thus, if $a\in
A_0$, then $a\bullet b=ab$ and $b\bullet a=ba$ for all $b\in A$.
Moreover, $\u{\s}(A)_0=A_0$ as algebras because $H^{\s}=H$ as
coalgebras and $\us(A)=A$ as right $H$-comodules.

\begin{lemma}\lelabel{3.14}
Let $A$ be a right $H^{op}$-comodule algebra. Then $\u{\s}(A)/\u{\s}(A)_0$ is
an $(H^{\s})^{op}$-Galois extension if and only if $A/A_0$ is an $H^{op}$-Galois extension.
\end{lemma}

\begin{proof} Note that $\u{\s}(A)\ot_{\u{\s}(A)_0}\u{\s}(A)=A\ot_{A_0}A$ as $k$-modules. We claim that the following diagram commutes:
$$\begin{array}{ccc}
\u{\s}(A)\ot_{\u{\s}(A)_0}\u{\s}(A)&\stackrel{\b^{\s}}{\longrightarrow}&
\u{\s}(A)\ot H^{\s}\\
\downarrow=&&\downarrow\xi\\
A\ot_{A_0}A&\stackrel{\b}{\longrightarrow}&A\ot H,
\end{array}$$
where $\b^{\s}$ is the Galois map for the right $(H^{\s})^{op}$-comodule algebra
$\u{\s}(A)$ and $\xi$ is the isomorphism given in \leref{3.13}.
Indeed, given $a$, $b\in\u{\s}(A)$, we have
$$\begin{array}{rcl}
(\xi\b^{\s})(a\ot b)&=&\xi(\S a\bullet b_{(0)}\ot b_{(1)})\\
&=&\S\xi(a_{(0)}b_{(0)}\ot b_{(2)})\s^{-1}(b_{(1)}\ot a_{(1)})\\
&=&\S(a_{(0)}b_{(0)}\ot b_{(6)})\s(S^{-1}(b_{(4)})\ot b_{(3)})\\
&&\s^{-1}(S^{-1}(b_{(5)})\ot b_{(1)}a_{(1)})\s^{-1}(b_{(2)}\ot a_{(2)})\\
&\stackrel{(\ref{rcocycle})}{=}&\S(a_{(0)}b_{(0)}\ot b_{(7)})\s(S^{-1}(b_{(4)})\ot b_{(3)})\\
&&\s^{-1}(S^{-1}(b_{(6)})b_{(1)}\ot a_{(1)})\s^{-1}(S^{-1}(b_{(5)})\ot b_{(2)})\\
&=&\S(a_{(0)}b_{(0)}\ot b_{(3)})\s^{-1}(S^{-1}(b_{(2)})b_{(1)}\ot a_{(1)})\\
&=&\S ab_{(0)}\ot b_{(1)}\\
&=&\b(a\ot b).\\
\end{array}$$
It follows that $\beta^{\s}$ is isomorphic if and only if $\beta$ is isomorphic. That is, $\us(A)/{\us(A)_0}$ is $({H^{\s}})^{op}$-Galois if and only if $A/A_0$ is $H^{op}$-Galois.
\end{proof}

\begin{theorem}\thlabel{3.15}
The following diagram commutes with exact rows:
$$\begin{array}{ccccccc}
1&\longrightarrow&{\rm Br}(k)&\longrightarrow&{\rm BC}(k, H, R)&
\stackrel{\stackrel{\sim}{\pi}}{\longrightarrow}&{\rm Gal}(\mathcal{H}_R)\\
&&\downarrow=&&\downarrow\u{\s}&&\downarrow\u{\s}\\
1&\longrightarrow&{\rm Br}(k)&\longrightarrow&{\rm BC}(k, H^{\s}, R^{\s})&
\stackrel{\stackrel{\sim}{\pi}}{\longrightarrow}&{\rm Gal}(\mathcal{H}^{\s}_{R^{\s}}).\\
\end{array}$$
\end{theorem}

\begin{proof}It is enough to verify that $\u{\s}\stackrel{\sim}{\pi}=\stackrel{\sim}{\pi}\u{\s}$
as group homomorphisms from ${\rm BC}(k,H,R)$ to ${\rm Gal}(\mathcal{H}^{\s}_{R^{\s}})$.
Let $A$ be a right $H^{op}$-comodule algebra. Suppose that $A$ is a Galois $R$-Azumaya algebra. By \coref{2.4} and \leref{3.14}
the right $(H^{\s})^{op}$-comodule algebra $\u{\s}(A)$ is a Galois
$R^{\s}$-Azumaya algebra. Since $\u{\s}(A)_0=A_0$ and $a\bullet b=ab$ for all $a\in A_0$ and $b\in \us(A)$, we have that $C_{\u{\s}(A)}(\u{\s}(A)_0)=C_A(A_0)$.
This means that $\pi(\u{\s}(A))=\pi(A)$ as right $H$-comodules since $H^{\s}=H$ as coalgebras. It follows that $\pi(\u{\s}(A))=\u{\s}(\pi(A))$ as right $(H^{\s})^{op}$-comodule algebras. It remains to verify  that the $H^{\s}$-actions on $\pi(\u{\s}(A))$ and $\u{\s}(\pi(A))$ coincide.

Let $\S X_i(h)\ot Y_i(h)=\b^{-1}(1\ot h)$ in $A\ot_{A_0}A$ for any $h\in H$. Then the $H$-action
on $\pi(A)$ is given by (\ref{MiU}). By \leref{3.13} and the proof of \leref{3.14}, we know that
$$(\b^{\s})^{-1}(1\ot h)=(\b^{-1}\xi)(1\ot h)=\S\s(S^{-1}(h_{(2)})\ot h_{(1)})X_i(h_{(3)})\ot Y_i(h_{(3)}),$$
for all $h\in H^{\s}$. Hence the $H^{\s}$-action on  $\pi(\u{\s}(A))$ is given by
$$h\cdot_{\s}a=\S\s(S^{-1}(h_{(2)})\ot h_{(1)})X_i(h_{(3)})\bullet a\bullet Y_i(h_{(3)}),
\hspace{0.2cm}h\in H^{\s}, a\in \pi(\u{\s}(A)).$$
If we define two right $H$-comodule structures on $A\ot_{A_0}A$ and $A\ot H^{op}$ by
$a\ot b\mapsto\S a\ot b_{(0)}\ot b_{(1)}$ and $a\ot h\mapsto\S a\ot h_{(1)}\ot h_{(2)}$
respectively, then it is easy to see that the Galois map $\b$ is an $H$-comodule isomorphism
with respect to the comodule structures. It follows that
$$\S X_i(h)\ot Y_i(h)_{(0)}\ot Y_i(h)_{(1)}=\S X_i(h_{(1)})\ot Y_i(h_{(1)})\ot h_{(2)}$$
for all $h\in H$. On the other hand,
if we define two right $H$-comodule structures on $A\ot_{A_0}A$ and $A\ot H^{op}$ by
$a\ot b\mapsto\S a_{(0)}\ot b\ot a_{(1)}$ and $a\ot h\mapsto\S a_{(0)}\ot h_{(2)}\ot S^{-1}(h_{(1)})a_{(1)}$
respectively, then  $\beta$ is also an $H$-comodule isomorphism
with respect to the foregoing defined comodule structures. It follows that
$$\S X_i(h)_{(0)}\ot Y_i(h)\ot X_i(h)_{(1)}=\S X_i(h_{(2)})\ot Y_i(h_{(2)})\ot S^{-1}(h_{(1)})$$
for all $h\in H$. Now we have
$$\begin{array}{rcl}
h\cdot_{\s}a&=&\S\s(S^{-1}(h_{(2)})\ot h_{(1)})X_i(h_{(3)})\bullet a\bullet Y_i(h_{(3)})\\
&=&\S\s(S^{-1}(h_{(2)})\ot h_{(1)})X_i(h_{(3)})\bullet(a_{(0)}Y_i(h_{(3)})_{(0)})
\s^{-1}(Y_i(h_{(3)})_{(1)}\ot a_{(1)})\\
&=&\S\s(S^{-1}(h_{(2)})\ot h_{(1)})X_i(h_{(3)})\bullet(a_{(0)}Y_i(h_{(3)}))
\s^{-1}(h_{(4)}\ot a_{(1)})\\
&=&\S\s(S^{-1}(h_{(2)})\ot h_{(1)})X_i(h_{(3)})_{(0)}a_{(0)}Y_i(h_{(3)})_{(0)}\\
&&\s^{-1}(Y_i(h_{(3)})_{(1)}a_{(1)}\ot X_i(h_{(3)})_{(1)})\s^{-1}(h_{(4)}\ot a_{(2)})\\
&=&\S X_i(h_{(4)})a_{(0)}Y_i(h_{(4)})\s(S^{-1}(h_{(2)})\ot h_{(1)})\\
&&\s^{-1}(h_{(5)}a_{(1)}\ot S^{-1}(h_{(3)}))\s^{-1}(h_{(6)}\ot a_{(2)})\\
&\stackrel{(\ref{+-})}{=}&\S X_i(h_{(5)})a_{(0)}Y_i(h_{(5)})
\s^{-1}(h_{(7)}a_{(2)}\ot S^{-1}(h_{(3)})h_{(2)})\\
&&\s(h_{(6)}a_{(1)}S^{-1}(h_{(4)})\ot h_{(1)})\s^{-1}(h_{(8)}\ot a_{(3)})\\
&=&\S X_i(h_{(3)})a_{(0)}Y_i(h_{(3)})
\s(h_{(4)}a_{(1)}S^{-1}(h_{(2)})\ot h_{(1)})\s^{-1}(h_{(5)}\ot a_{(2)})\\
&\stackrel{(\ref{MiU})}{=}&\S h_{(3)}\cdot a_{(0)}
\s(h_{(4)}a_{(1)}S^{-1}(h_{(2)})\ot h_{(1)})\s^{-1}(h_{(5)}\ot a_{(2)})\\
&\stackrel{(\ref{sigmaaction})}{=} & h\rhu a.
\end{array}$$
Thus we have proved that  $h\cdot_{\s}a=h\rhu a$ for all $a\in\pi(\u{\s}(A))=\u{\s}(\pi(A))$
and $h\in H^{\s}$. Whence, $\pi(\u{\s}(A))=\u{\s}(\pi(A))$ as {\rm YD} $H^{\s}$-module algebras.
\end{proof}

\section*{ACKNOWLEDGMENTS}
\hskip\parindent
The first author would like to thank the School of Mathematics, Statistics and Computer Science, Victoria University of Wellington for their hospitality during his visit in 2004. He is grateful to the URF of VUW for the financial support. He is also supported by NSF of China (No. 10471121).\\
The second author is supported by the Marsden Fund.

\end{document}